\documentclass{amsart}

\usepackage{graphicx}
\usepackage[colorlinks=true, allcolors=blue]{hyperref}

\usepackage{comment}

\usepackage[letterpaper,top=2cm,bottom=2cm,left=3cm,right=3cm,marginparwidth=1.75cm]{geometry}
\usepackage{amsmath}

\usepackage{amsthm}

\usepackage{amsfonts,a4wide}
\usepackage{mathrsfs}
\usepackage{amssymb}
\usepackage{mathtools}
\usepackage{tikz}
\usepackage{rotating}
\usetikzlibrary{shapes}
\newtheorem{theorem}{Theorem}
\newtheorem{corollary}[theorem]{Corollary}
\newtheorem{lemma}[theorem]{Lemma}
\newtheorem{proposition}[theorem]{Proposition}

\newtheorem{definition}[theorem]{Definition}
\newtheorem{remark}[theorem]{Remark}
\newtheorem{example}[theorem]{Example}

\newcommand{\spa}{\mbox{\rm span}}

\newcommand{\Real}{\mathbb{R}}

\newcommand{\Zzz}{\mathbb{Z}}
\newcommand{\T}{\mathbb{T}}

\newcommand{\h}{\mathcal{H}}

\newcommand{\seq}[1]{\left<#1\right>}

\usepackage{amsmath}
\usepackage{graphicx}
\usepackage[colorlinks=true, allcolors=blue]{hyperref}
\usepackage{amsfonts}

\usepackage{color}

\title
{An operator theory approach to the evanescent part of a two-parametric weak-stationary stochastic process}
\author{Zbigniew Burdak, Marek Kosiek, Patryk Pagacz, Marek S\l oci\'nski
}

\begin{document}

\begin{abstract}
   A new approach to the evanescent part of a two-dimensional  weak-stationary stochastic process with the past given by a half-plane is proceed.
   The classical result due to Helson and Lowdenslager divides a two-parametric weak-stationary stochastic process into three parts. In this paper we describe the most untouchable one - the evanescent part. Moreover, we point out   how this part depends on the shape of the past.
\end{abstract}
\subjclass{MSC2020: 47B20, 47A13,60G10,60G25}%

\keywords{Keywords: evanescent part, stochastic process, weak-stationary, isometies, compatible isometries}

\maketitle

\section{Introduction}
 In the remarkable work \cite{W} 
 Wold decomposed a stationary stochastic process between a deterministic process and the moving average of a white noise. A stochastic process has an interpretation as a Hilbert space isometry. Halmos gave an abstract version of Wold's result in \cite{Halmos}   decomposing an arbitrary Hilbert space isometry between a unitary operator and a unilateral shift with a geometrical description of a unilateral shift. The similar decomposition of an isometry without the geometrical structure was pointed out earlier by  
von Neumann in \cite{vN}.

Another direction of generalization of Wold's result are multi-dimensional stochastic processes called random fields. Random fields find their practical meaning, in the context presented in this paper, e.g in the problem of image texture modeling, filtering data and other applications, see \cite{RANradar,C,app,F,FNradar,FMP,FNW}. The set of indices of a stochastic process is usually interpreted as time and is divided among a past, a present and a future. In the one-dimensional case negative indices naturally define the past.

In the two-dimensional case there are at least two concepts of the definition of the past: a quaternary past and a half-plane past. Moreover, a random field strongly depends on the definition of the past.  Helson and Lowdenslager in \cite{HL1961} decomposed an arbitrary random field with a half-plane past among three summands, where a precise description of so called evanescent part is still missing. An abstract version of this result for commuting pairs of Hilbert space isometries is given in \cite{S}.

The paper is organized as follows. In Subsection \ref{PT} we recall the basic concepts 
of prediction theory introduced in \cite{HL,HL1961}. In Subsections \ref{OT} we present the operator theory approach to weak-stationary random fields as pairs of commuting Hilbert space isometries. In particular we include the recent  results on decomposition of pairs of commuting Hilbert space isometries and models of summands in the decomposition. In Sections \ref{SectL}, \ref{SectRS} and \ref{SectNS} we  show that the evanescent part of a random field with a half-plane past corresponds to exactly one pair of isometries (one summand) in the decomposition, where the corresponding summand is determined by a half-plane defining the past as follows: horizontal half-plane - a pair of a unitary operator and a unilateral shift, a half-plane with rational slope -  a pair of generalized powers, and a half-plane with irrational slope - a coninuously given pair.

Additionally, in the appendix, we give a formal proof of the form of a half-plane. Such result was mentioned in \cite{Cuny}, but we were not able to find a strict proof in  literature.

\subsection{Prediction theory}\label{PT}

A random field is a collection of random variables $\{X_{t}\}_{t\in T}$ defined on the same probability space indexed in some topological space $T$ usually interpreted as time. We consider a complex-valued variables indexed by $T=\Zzz^2$.
\begin{definition}\label{weakstat}
  A random field $\{X_{(s,t)}\}_{(s,t)\in \Zzz^2}$ is \emph{weak-stationary} if:
\begin{itemize}
\item $\mathbb{E}[|X_{(s,t)}|^2]<\infty$ for any $(s,t)\in \Zzz^2$,
\item $\mathbb{E}[X_{(s,t)}]$ do not depend on $(s,t)\in \Zzz^2$ (frequently assumed to vanish),
\item the cross moments depend at most on the distance i.e. $\mathbb{E}[X_{(s,t)}\overline{X_{(s,t)+(p,q)}}]=\gamma(p,q)$, for some function $\gamma$.
\end{itemize}
\end{definition}
By the first condition a weak-stationary random field can be considered in a Hilbert space of square integrable functions with an inner product $\seq{X,Y}=\mathbb{E}[X\overline{Y}]$.

Let $\{X_{(s,t)}\}_{(s,t)\in \Zzz^2}$ be a weak-stationary random field. The main goal of prediction theory is to investigate its dependence on the past. In other words,
dependence of $X_{(s,t)}$ on $\{X_{(i,j)}: (i-s,j-t)\in S\}$ where the subset $S\subset \Zzz^2$ define "the past". Let us recall precisely two concepts of a two-dimensional past:
\begin{equation}\label{QL}
Q =\{(i,j) : i\leq 0,j\leq 0\}\setminus \{(0,0)\},\quad  %L =\{(i,j) : i\leq -1,j\in\Zzz\}\cup\{(0,j) :j\leq -1\},
L =\{(i,j) : j\leq -1,i\in\Zzz\}\cup\{(i,0) :i\leq -1\}
\end{equation}
which correspond to the third-quadrant and the bottom half-plane, respectively.
The set $L$ is a special case of a half-plane as defined in \cite{HL}.

\begin{definition}\label{hp}
$S$ is a \emph{half-plane} of lattice points if
\begin{itemize}
\item $(0,0)\notin S$,
\item $S$ is a semigroup,
\item$(m,n)\in S$ if and only if $(-m, -n)\not\in S $ unless $m=n=0$.
\end{itemize}
\end{definition}

There is one-to-one correspondence between half-planes and (preserving addition) total orders on $\Zzz^2$. Indeed, a half-plane $S$ define a total order as follows: $$(s,t)\prec (i,j) \Longleftrightarrow (s,t)\in (i,j)+S\cup \{(0,0)\}.$$  Conversely,   a half-plane can be defined by a total order $\prec$ as follows: $$S:=\{(i,j)\in \Zzz : (i,j)\prec(0,0), (i,j)\not=(0,0)\}.$$

In \cite{Cuny} there was observed that a half-plane can be described by a vector $\textbf{v}$, i.e. each half-plane is equal to 
\begin{equation}\label{S_v}
S_{\textbf{v}} =\{(i,j) : \seq{\textbf{v},(i,j)}=iv_1+jv_2\geq 0 \textnormal{ for } i<0, \seq{\textbf{v},(i,j)} > 0 \textnormal{ for } i\geq 0\},
\end{equation} or to \begin{equation}\label{S_vhat}\widehat{S}_{\textbf{v}} =\{(i,j) : \seq{\textbf{v},(i,j)}=iv_1+jv_2> 0 \textnormal{ for } i\le0, \seq{\textbf{v},(i,j)} \ge 0 \textnormal{ for } i> 0\}.\end{equation} for some $\textbf{v}\in \Real^2\setminus \{(0,0)\}$. It can be showed that $\widehat{S}_{[x,y]}=\{(i,j): (j,i)\in S_{[y,x]}\}.$
The formal prove that any half-plane can be described as above is posted in Appendix. In particular, set $L$ is equal to $S_{[0,-1]}$. Thus half-planes are rotations of $L$. 
The half-planes $S_{\textbf{v}}$ were also called a \emph{nonsymmetrical or augmented half-planes} f.e. in \cite{ChiangJMA,ChiangAMAS,F}.

The general prediction theory based on weak-stationary random fields with a half-plane past was developed in the seminal papers by Helson and Lowdenslager (see \cite{HL,HL1961}). Let us recall some results following from their work.

The function $\gamma$ on $\Zzz^2$ corresponding to a weak-stationary random field in Definition \ref{weakstat} is  positive-definite. Hence, by Bochner's theorem, there is a unique distribution function $F$ on the torus $\T=(-\pi,\pi]^2$ such that
$$\gamma(s,t) =\int_{-\pi}^{\pi}\int_{-\pi}^{\pi} exp(-i(s\lambda_1 +t\lambda_2))dF(\lambda_1,\lambda_2).$$ 
Further, let $H:=\spa{\{X(s,t); (s,t) \in\Zzz^2\}}$ where $\spa\{\dots\}$ denotes a closed linear manifold and $L^2(dF)$ be the space of square integrable functions on $\T$ with respect to the measure $dF(\lambda_1,\lambda_2)$. Both spaces $H, L^2(dF)$ are Hilbert spaces and the map $X_{(s,t)} \to exp(-i(s\lambda_1 +t\lambda_2))$ extends to an isomorphism from $H$ onto $L^2(dF)$.

 Let \begin{equation}\label{Hij}H^{(i,j)}:=\spa\{X_{(s,t)}:(s-i,t-j) \in S\}\end{equation} be the past of $X_{(i,j)}$. In the context of Helson and Lowdenslager paper \cite{HL1961} we give a standard and meaningful definition.
\begin{definition}\label{processesDef}
A random field $\{X_{(s,t)}\}_{(s,t)\in \Zzz^2}$ is said to be
\begin{itemize}
  \item \emph{deterministic} if $X_{(s,t)}\in \bigcap\limits_{(i,j)\in\Zzz^2}H^{(i,j)}$, for $(s,t)\in\Zzz^2$,
  \item \emph{evanescent} if $ \bigcap\limits_{(i,j)\in\Zzz^2}H^{(i,j)}=\{0\}$ and $X_{(s,t)}\in H^{(s,t)}$ for $(s,t)\in \Zzz^2$,
  \item \emph{purely nondeterministic} if $X_{(s,t)} \in \spa\{I_{(i,j)}:(i,j)\in (s,t)+S\cup\{(0,0)\}\}$, for $(s,t)\in\Zzz^2$, where
   $I_{(i,j)}:=X_{(i,j)}-P_{H^{(i,j)}}X_{(i,j)}$ is called \emph{an innovation part} of $X_{(i,j)}$, for $(i,j)\in\Zzz^2$.
\end{itemize}

\end{definition}
The decomposition among the three parts above is given in \cite{HL1961}. The version we present follows from \cite[Theorem 5]{HL1961} and \cite[Theorems 3.1, 3.2]{Cuny} (in here for scalar valued process). 

\begin{theorem}\label{WoldbyCuny}
Let $\{X_{(s,t)}\}_{(s,t)\in\Zzz^2}$ be a weak-stationary random field with a half-plane past $S$ and a unique distribution function $F$ with Lebesgue decomposition $F=fdm_2 +dF_s$. The process $\{X_{(s,t)}\}_{(s,t)\in\Zzz^2}$ has a nontrivial innovation part if and only if $log(f)\in L^1(dm_2)$ and then:
\begin{itemize}
\item
\begin{equation}\label{movavr}
    A_{(s,t)}=\sum_{(k,l)\in S\cup\{(0,0)\}} a_{k,l}I_{s+k,t+l},
  \end{equation}
  is a purely nondeterministic random field with distribution  $dF_A=fdm_2$ where $I_{(s,t)}$ is an innovation part of $X_{(s,t)}$ and $a_{k,l}=\seq{I_{0,0},I_{0,0}}^{-1}\seq{X_{0,0},I_{k,l}}$, for $(k,l)\in S\cup \{(0,0)\}$,
\item
 there are an evanescent random field $\{E_{(s,t)}\}_{(s,t)\in\Zzz^2}$ with distribution $F_E$ and a deterministic random field $\{D_{(s,t)}\}_{(s,t)\in\Zzz^2}$ with distribution $F_D$ such that $F_s=F_E+F_D$ and
$$ X_{(s,t)}=A_{(s,t)}+E_{(s,t)}+D_{(s,t)}, \quad (s,t)\in\Zzz^2.$$
\end{itemize}

 Moreover, there is a Borel set $\Gamma\subset \T$  of Lebesgue measure $0$ such that

  $$A_{(s,t)}=\int_{\T\setminus \Gamma}\exp(is\lambda_1+it\lambda_2)dF X_{(0,0)}.$$
\end{theorem}

The main difference between the above two-dimensional Wold-type decomposition and the classical Wold's one-dimensional version is the evanescent part. Let us point out that a purely nondeterministic part does not depend on the choice of $S$ but evanescent and deterministic parts do. Recent spectral approaches to the evanescent part may be found in \cite{C,CF,Cuny,FMP,KMN,KF,KL,L,MSS}.

\subsection{Two-dimensional Wold decomposition in operator theory}\label{OT}

There is a natural interpretation of a random field $\{X_{(s,t)}\}_{(s,t)\in \Zzz^2}$ in the Hilbert space of square integrable functions with the inner product $\seq{X,Y}:= \mathbb{E}(X\overline{Y})$ as $H=\spa\{X_{(i,j)}:(i,j)\in\Zzz^2\}$. Then operators defined by
\begin{equation}\label{processToOp}
U_1X_{(s,t)}=X_{(s-1,t)}, \quad U_2X_{(s,t)}=X_{(s,t-1)}, \quad (s,t)\in \Zzz^2
\end{equation}
are unitary on $H$ and commute. Since without lost of generality we may assume that $\mathbb{Z}_-^2\subset S$ subspaces $H^{(i,j)}$ defined in \eqref{Hij} are invariant under $U_1,U_2,$ so \begin{equation}\label{processToOpV}V_1=U_1|_{H^{(i,j)}}, V_2=U_2|_{H^{(i,j)}}\end{equation} are commuting isometries, not necessarily unitary for any $(i,j)\in\mathbb{Z}^2$. Operators $V_1,V_2$ clearly depend on $(i,j)$ but usually $(i,j)$ is clear from the context and we skip it in notation.

As one can expect a decomposition of a random field provides a decomposition of $(V_1, V_2)$. It turns out that there is also the reverse connection. Precisely random fields may be investigated as pairs of commuting isometries. Let us describe this relations. The starting point is an operator version of Helson and Lowdenslager decomposition obtained by I. Suciu in \cite{S}. 

Recall the operator theory notions corresponding to deterministic, evanescent and purely non-deterministic random fields.  A pair of isometries (operators) is \emph{completely non-unitary} if there is no non-trivial subspace reducing both operators to unitary operators. Further:

\begin{definition}\label{defSuciu}
Let $S\subset\Zzz^2$ be a semigroup such that $S\cap -S =\emptyset$ and $\Zzz_+^2\subset S$ and $(V_1,V_2)$ be a pair of commuting isometries on a Hilbert space $H$ such that $$V^{(m,n)}:=\begin{cases}
                       V_1^mV^n_2, & \mbox{if } m,n \geq 0,\\
                       V_1^{*|m|}V^n_2, & \mbox{if } m<0\leq n, \\
                       V_2^{*|n|}V^m_1, & \mbox{if } n<0\leq m,                     \end{cases}$$
 is an isometry for any $(m,n)\in S$.

Then $(V_1,V_2)$ is called
\begin{itemize}
  \item \emph{quasi-unitary} if $$cl\left\{\bigcup_{(m,n)-(k,l)\in S} V^{(k,l)*}V^{(m,n)}H\right\}=H,$$
  \item \emph{totally non-unitary} if there is no non-trivial subspace reducing $(V_1,V_2)$ to a quasi-unitary pair,
  \item \emph{strange} if it is quasi-unitary and completely non-unitary.
\end{itemize}
\end{definition}

Note that $\{V^{(m,n)}\}_{(m,n)\in S}$ is a semi-group since $V^{(m,n)}$ is assummed to be an isometry for any $(m,n)\in S.$

\begin{theorem}\label{SuciuMain}{\cite[Theorem 3]{S}}
 Let $S\subset\Zzz^2$ be a semigroup such that $S\cap -S =\emptyset$ and $\Zzz_+^2\subset S$. Let a pair of isometries $(V_1,V_2)$ on a Hilbert space $H$ be such that $V^{(m,n)}$ is an isometry for any $(m,n)\in S$.
 Then there is a unique decomposition
 $$ H=H_u\oplus H_e \oplus H_t,$$
 where  $H_u, H_e, H_t$ are $(V_1,V_2)$ reducing subspaces, $(V_1|_{H_u},V_2|_{H_u})$ is a pair of unitary operators, $(V_1|_{H_e},V_2|_{H_e})$ is strange and $(V_1|_{H_t},V_2|_{H_t})$ is totally non-unitary.
\end{theorem}
Theorem \ref{SuciuMain} implies Theorem \ref{WoldbyCuny} by Remark \ref{Back} below. First we remark some property used in the proof of Remark \ref{Back} and further in the paper.
\begin{remark}\label{projprop}
Let $(V_1, V_2)$ be a pair of isometries on $H$ and $(U_1, U_2)$ be its minimal unitary extension on $K$. If $H=H_1\oplus H_2$ is a decomposition between $(V_1, V_2)$ reducing subspaces, then $K=K_1\oplus K_2$ where $K_i$ are the spaces of the minimal unitary extensions of $(V_1|_{H_i}, V_2|_{H_i})$ for $i=1,2$. It follows a decomposition $K=H_1\oplus (K_1\ominus H_1)\oplus H_2\oplus (K_2\ominus H_2)$ which is equivalent to $P_{H_i}=P_{H}P_{K_i}=P_{K_i}P_{H}$ for $i=1,2$ as operators on $K.$ In other words $H_i=K_i\cap H$ for $i=1,2$.

Clearly the similar result is correct for decomposition among more reducing subspaces.
\end{remark}

\begin{remark}\label{Back}
Let us consider decompositions $$H^{(i,j)}=H_u^{(i,j)} \oplus H_e^{(i,j)}\oplus H_t^{(i,j)}$$ of subspaces \eqref{Hij} given by Theorem \ref{SuciuMain} for $(V_1,V_2)$  as in  \eqref{processToOpV} for $(i,j)\in\mathbb{Z}^2$. A minimal unitary extension of $(V_1,V_2)$ for each $(i,j)$ is an extension to the whole $H$ and by Remark \ref{projprop} we get $H=H_u\oplus H_e\oplus H_t$ where $H_u, H_e, H_t$ are minimal unitary extensions of $H_u^{(i,j)}, H_s^{(i,j)}, H_t^{(i,j)}$ respectively. Clearly $H_u=H_u^{(i,j)}$ so it does not depend on $(i,j)$. Moreover, since by Remark \ref{projprop} $P_{H_t^{(i,j)}}=P_{H_t}P_{H^{(i,j)}}$ and $P_{H^{(i,j)}}\xrightarrow{SOT} I$ for $(i,j)\rightarrow (\infty,\infty)$ we get $H_t=\spa\{H_t^{(i,j)}:(i,j)\in\mathbb{Z}^2\}$ and similarly $H_e= \spa\{H_e^{(i,j)}:(i,j)\in\mathbb{Z}^2\}.$

Moreover, if $$X_{(s,t)}=A_{(s,t)}+E_{(s,t)}+D_{(s,t)}$$ is a decomposition like in Theorem \ref{WoldbyCuny}, then $D_{(s,t)}\in H_u,\; A_{(s,t)}\in H_t,\;E_{(s,t)}\in H_e$. Hence Theorem \ref{SuciuMain} implies Theorem \ref{WoldbyCuny} by $$X_{(s,t)}=P_{H_t}X_{(s,t)}+P_{H_e}X_{(s,t)}+P_{H_u}X_{(s,t)}.$$
\end{remark}

The Wold(-type) decompositions of Hilbert space isometries are originated as generalization of the results for stochastic processes/random fields. Hence the result as in Remark \ref{Back} is not surprising. However, operator theory approach to the problem developed independently. In this paper we apply some recent result in operator theory concerning pairs of commuting, compatible isometries, to get a precise description of the evanescent part in the case of \emph{rational nonsymmetrical half-plane} (RNHP) i.e. $S=S_{\textbf{v}}$, where $\textbf{v}=(k,l)\in \Zzz^2$ and also some result describing the evanescent part for irrational  nonsymmetrical half-plane past i.e. $S=S_{\textbf{v}}$, where $\textbf{v}=(q,r)\in \mathbb{R}^2,\; \frac{q}{r}\notin\mathbb{Q}.$ The problem of describing the evanescent part under above assumption on the past was investigated by many authors (i.e. \cite{RANradar,Cuny,F,FMP,FNradar,FNW,KMN,KF,KL,L}). 

The pair of isometries $(V_1, V_2)$ on $\mathcal{H}$ is compatible if $P_{V_1^m(\mathcal{H})}P_{V_2^n(\mathcal{H})}=P_{V_2^n(\mathcal{H})}P_{V_1^m(\mathcal{H})}$ for any non negative $m,n$ (\cite{Mil}). Let $(V_1, V_2)$ be as in \eqref{processToOpV} for any $(i,j)\in\mathbb{Z}^2$. Note that $V_1^m(H^{(i,j)})=U_1^m(H^{(i,j)})=H^{(i-m,j)}$ and $V_2^m(H^{(i,j)})=U_2^n(H^{(i,j)})=H^{(i,j-n)}$. On the other hand, a half-plane past define a total order which implies a total inclusion order among subspaces $H^{(i,j)}$ and in turn a total order among $P_{H^{(i,j)}}$. It clearly implies compatibility of  $(V_1, V_2)$. Let us remark that a pair of isometries corresponding to random field with $Q$ past \eqref{QL} does not have to be compatible.

The commuting, compatible pairs are described in \cite{BKS, BKPSCorr}. More precisely, from \cite[Theorem 3]{BKPSCorr} we get a decomposition among pairs of isometries of certain types. First we formulate a decomposition theorem. Then we recall definitions of the types of compatible pairs that appears in the decomposition and its relation with the random fields.
\begin{theorem}\label{dec_thm}
For any pair of commuting, compatible isometries $(V_1, V_2)$ on  the Hilbert space $H$ there is a decomposition:
$$H=H_{uu}\oplus H_{us}\oplus H_{su}\oplus H_{gp}\oplus H_d\oplus H_c$$
where $H_{uu}, H_{us}, H_{su}, H_{gp}, H_d, H_c$ reduce $(V_1, V_2)$ such that:
\begin{itemize}
\item $V_1|_{H_{uu}}, V_2|_{H_{uu}}, V_1|_{H_{us}}, V_2|_{H_{su}}$ are unitary and $V_1|_{H_{s u}}, V_2|_{H_{u s}}$ are unilateral shifts,
\item $(V_1|_{H_{gp}}, V_2|_{H_{gp}})$ decomposes among pairs of generalized powers,
\item $(V_1|_{H_d}, V_2|_{H_d})$ decomposes among pairs given by diagrams,
\item $(V_1|_{H_c}, V_2|_{H_c})$ is a continuously given pair.
\end{itemize}
\end{theorem}

The deterministic part is described by Theorem \ref{SuciuMain} as a pair of unitary operators.

The purely nondeterministic part is described as a pair defined by a diagram by Theorem \ref{SuciuTw4}.
A \emph{diagram} is a set  $J\subset \Zzz^2$ such that $J+\Zzz^2_+\subset J$. The concept of  pairs of isometries defined by diagrams was introduced by M\" uller in \cite{Mil} and developed in \cite{BKPS}.

\begin{definition}
\label{on_Diagram}
Let $J$ be a diagram and $\mathcal{H}$ be a Hilbert space. Define $$H_J:=\{f\in L^2_{\mathcal{H}}(\T): \hat{f}_{i,j}=0\text{ for }(i,j)\notin J\}=\spa\{z_1^iz_2^jh:(i,j)\in J, h\in\h\}$$ where  $L^2_{\mathcal{H}}(\T)$ denotes a space of $\mathcal{H}$ valued, square integrable functions and $\hat{f}_{i,j}$ denotes the respective Fourier coefficient. Note that $H_J$ is invariant under operators of multiplication by independent variables $M_{z_1}, M_{z_2}$.

A pair of \emph{isometries defined by the diagram $J$ and the space
$\mathcal{H}$} is a pair unitarily equivalent to $M_{z_1}|_{H_J}, M_{z_2}|_{H_J}.$
\end{definition}

For the one who prefers purely geometrical description let us recall an equivalent definition as in \cite{BKS}.

\begin{remark}\label{diagram_equiv}
Let $J$ be a diagram and $\mathcal{H}$ be a Hilbert space.

Define a Hilbert space $H:=\{\mathbf{v}=\{v_{(i,j)}\}_{(i,j)\in J}: v_{i,j}\in\mathcal{H}\}$ and operators

$V_{1}\mathbf{v}=\mathbf{x}$ where $x_{(i,j)}=\left\{
                                              \begin{array}{ll}
                                                v_{(i-1,j)}, & \hbox{if $(i-1,j)\in J$;} \\
                                                0, & \hbox{$(i-1,j)\notin J$.}
                                              \end{array}
                                            \right.$

$V_{2}\mathbf{v}=\mathbf{y}$ where $y_{(i,j)}=\left\{
                                              \begin{array}{ll}
                                                v_{(i,j-1)}, & \hbox{if $(i,j-1)\in J$;} \\
                                                0, & \hbox{if $(i,j-1)\notin J$.}
                                              \end{array}\right.$

Operators $V_1, V_2$ are commuting isometries defined by the diagram $J$ and the space $\mathcal{H}.$

\end{remark}
In particular for $J=(\Zzz_+\cup\{0\})^2$ the subspace $H_J$ is a (vector valued) Hardy space and $M_{z_1}|_{H_J}, M_{z_2}|_{H_J}$ is a pair of doubly commuting unilateral shifts. Let us note that the set $(i,j)+J$, called a \emph{translation} of the diagram $J$, is a diagram which generate the same pair of isometries as $J$ (assuming the same space $\mathcal{H}$),  where $(i,j)\in\Zzz^2.$ Hence, it is more precise to view a pair of isometries as defined by a translation equivalence class of a diagram than a single diagram.

Since $S$ in Definition \ref{defSuciu} is a semigroup containing $\Zzz_+^2$ it is a diagram. Hence, directly from \cite[Theorem 4]{S} follows.

\begin{theorem}\label{SuciuTw4}
  Any totally non-unitary pair of isometries is a pair given by a diagram.
\end{theorem}
Note that Theorem \ref{SuciuTw4} implies \eqref{movavr}.

 The evanescent part of the random field with the past $L$ \eqref{QL}  is described by a pair of a unitary operator and a unilateral shift by Theorem \ref{Thm_HtHe}. The evanescent part of the random field with the rational nonsymmetrical half-plane past $S_{(k,l)}$ for $(k,l)\in\mathbb{Z}^2$  is described by a pair of generalized powers by Theorem \ref{Thmrs}. Generalized powers are generalization of a pair $(V^m, V^n)$ where $V$ is a unilateral shift and $m,n$ are positive integers. Hence a random field corresponding to $(V^m, V^n)$ is defined by a one-dimensional white noise. The model of a pair of generalized powers and its definition require the concept of a periodic diagram.

\begin{definition}\label{periodic diagram definition}
 The diagram $J$ is \emph{periodic} if there exist  $m,n$ such that for
 $J_0:=(\{0,1,\dots,m-1\}\times \mathbb{Z})\cap J$ and $J_k:=(km,-kn)+J_0=\{(i+km, j-kn):(i,j)\in J_0\}$ for $k \in \mathbb{Z}$  we have
 $J=\sum_{k\in\mathbb Z} J_k$,  and $J_k$ are disjoint for different $k$.
 Then the set $J_0$ is called a \emph{period} of the diagram.
 \end{definition}
 The idea of a periodic diagram and its period is explained in the following picture:

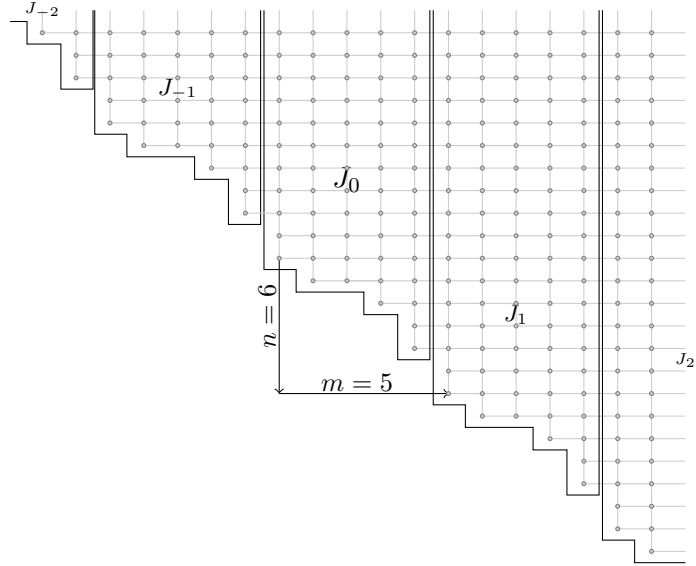
\begin{center}
 \begin{tikzpicture}
[yscale=0.3, xscale=0.45, auto,
kropka/.style={draw=black!50, circle, inner sep=0pt, minimum size=1.5pt, fill=black!20}]

\draw [black!20] (6,-11) -- (5,-11);
\draw [black!20] (6,-10) -- (4,-10);
\draw [black!20] (6,-9) -- (4,-9);
\draw [black!20] (6,-8) -- (3,-8);
\draw [black!20] (6,-7) -- (3,-7);
\draw [black!20] (6,-6) -- (3,-6);
\draw [black!20] (6,-5) -- (2,-5);
\draw [black!20] (6,-4) -- (0,-4);
\draw [black!20] (6,-3) -- (-1,-3);
\draw [black!20] (6,-2) -- (-2,-2);
\draw [black!20] (6,-1) -- (-2,-1);
\draw [black!20] (6,0) -- (-2,0);
\draw [black!20] (6,1) -- (-4,1);
\draw [black!20] (6,2) -- (-5,2);
\draw [black!20] (6,3) -- (-6,3);
\draw [black!20] (6,4) -- (-6,4);
\draw [black!20] (6,5) -- (-7,5);
\draw [black!20] (6,6) -- (-7,6);
\draw [black!20] (6,7) -- (-9,7);
\draw [black!20] (6,8) -- (-10,8);
\draw [black!20] (6,9) -- (-11,9);
\draw [black!20] (6,10) -- (-12,10);
\draw [black!20] (6,11) -- (-12,11);
\draw [black!20] (6,12) -- (-13,12);

\draw [black!20] (5,-11) -- (5,13);
\draw [black!20] (4,-10) -- (4,13);
\draw [black!20] (3,-8) -- (3,13);
\draw [black!20] (2,-5) -- (2,13);
\draw [black!20] (1,-5) -- (1,13);
\draw [black!20] (0,-4) -- (0,13);
\draw [black!20] (-1,-3) -- (-1,13);
\draw [black!20] (-2,0) -- (-2,13);
\draw [black!20] (-3,1) -- (-3,13);
\draw [black!20] (-4,1) -- (-4,13);
\draw [black!20] (-5,2) -- (-5,13);
\draw [black!20] (-6,4) -- (-6,13);
\draw [black!20] (-7,6) -- (-7,13);
\draw [black!20] (-8,7) -- (-8,13);
\draw [black!20] (-9,7) -- (-9,13);
\draw [black!20] (-10,8) -- (-10,13);
\draw [black!20] (-11,9) -- (-11,13);
\draw [black!20] (-12,10) -- (-12,13);
\draw [black!20] (-13,12) -- (-13,13);

\draw [->,black!100] (-6,2) -- (-6,-4);
\draw [->,black!100] (-6,-4) -- (-1,-4);
\node at (-6.1,-2) {\begin{rotate}{90}$n=6$\end{rotate}};
\node at (-3.7,-3.5) {$m=5$};

\draw [black!100] (-6.45,13) --  ++(0,-11.5) -- ++ (0.95,0)-- ++(0,-1) -- ++(2,0) -- ++(0,-1) -- ++(1,0)-- ++(0,-2) -- ++(0.95,0) -- ++(0,15.5);
\draw [black!100] (-11.45,13) --  ++(0,-5.5) -- ++ (0.95,0)-- ++(0,-1) -- ++(2,0) -- ++(0,-1) -- ++(1,0)-- ++(0,-2) -- ++(0.95,0) -- ++(0,9.5);
\draw [black!100] (-1.45,13) --  ++(0,-17.5) -- ++ (0.95,0)-- ++(0,-1) -- ++(2,0) -- ++(0,-1) -- ++(1,0)-- ++(0,-2) -- ++(0.95,0) -- ++(0,21.5);
\draw [black!100] (-13.95,12.5)  -- ++(0.5,0) -- ++(0,-1) -- ++(1,0)-- ++(0,-2) -- ++(0.95,0) -- ++(0,3.5);
\draw [black!100] (3.55,13) --  ++(0,-23.5) -- ++ (0.95,0)-- ++(0,-1) -- ++(1.5,0) ;

\draw [black!20] (-1,-3) --++(0,-1) -- ++(1,0) -- ++(0,-1) -- ++ (2,0) -- ++(0,-1) -- ++ (1,0)-- ++(0,-1) ;
\draw [black!20] (-2,-2) --  ++(0,2) -- ++ (-1,0)-- ++(0,1) -- ++(-2,0) -- ++(0,1) -- ++(-1,0) -- ++(0,1) ;
\draw [black!20] (-6,3)  --  ++(0,1) -- ++ (-1,0)-- ++(0,2) -- ++(-1,0) -- ++(0,1) -- ++(-2,0) -- ++(0,1) -- ++(-1,0) -- ++(0,1) ;
{\large \path (-4,5.5) node{$J_0$};}
{\small \path (-9,9.5) node{$J_{-1}$};}
{\small \path (1,-0.5) node{$J_1$};}
{\tiny \path (-13,13) node{$J_{-2}$};}
{\tiny \path (6,-2.5) node{$J_2$};}

\path (5,-11)  node[kropka]{} ++(0,1) node [kropka] {} ++(0,1) node [kropka] {}++(0,1) node [kropka] {}++(0,1) node [kropka] {}++(0,1) node [kropka] {}++(0,1) node [kropka] {}++(0,1) node [kropka] {}++(0,1) node [kropka] {}++(0,1) node [kropka] {}++(0,1) node [kropka] {}++(0,1) node [kropka] {}++(0,1) node [kropka] {}++(0,1) node [kropka] {}++(0,1) node [kropka] {}++(0,1) node [kropka] {}++(0,1) node [kropka] {}++(0,1) node [kropka] {}++(0,1) node [kropka] {}++(0,1) node [kropka] {}++(0,1) node [kropka] {}++(0,1) node [kropka] {}++(0,1) node [kropka] {}++(0,1) node [kropka] {};

\path (4,-10)  node[kropka]{} ++(0,1) node [kropka] {} ++(0,1) node [kropka] {}++(0,1) node [kropka] {}++(0,1) node [kropka] {}++(0,1) node [kropka] {}++(0,1) node [kropka] {}++(0,1) node [kropka] {}++(0,1) node [kropka] {}++(0,1) node [kropka] {}++(0,1) node [kropka] {}++(0,1) node [kropka] {}++(0,1) node [kropka] {}++(0,1) node [kropka] {}++(0,1) node [kropka] {}++(0,1) node [kropka] {}++(0,1) node [kropka] {}++(0,1) node [kropka] {}++(0,1) node [kropka] {}++(0,1) node [kropka] {}++(0,1) node [kropka] {}++(0,1) node [kropka] {}++(0,1) node [kropka] {};

\path (-2,-2)  node[kropka]{} ;

\path (-2,-1)   node [kropka]{}  ;

\path (-2,0)  node [kropka]{}   ++(-1,0) node [kropka] {}  ;

\path (-2,1) node [kropka] {}   ++(-1,0) node [kropka] {}     ++(-1,0) node [kropka] {}  ++(-1,0) node [kropka] {}  ;

\path (-2,2) node [kropka] {}   ++(-1,0) node [kropka] {}    ++(-1,0) node [kropka] {}     ++(-1,0) node [kropka] {}   ++(-1,0) node [kropka] {} ;

\path (-2,3)  node [kropka] {}   ++(-1,0) node [kropka] {}     ++(-1,0) node [kropka] {}     ++(-1,0) node [kropka] {}    ++(-1,0) node [kropka] {} ;

\path (-2,4)  node [kropka] {}  ++(-1,0) node [kropka] {}     ++(-1,0) node [kropka] {}     ++(-1,0) node [kropka] {}   ++(-1,0) node [kropka] {} ;

\path (-2,5) node [kropka] {}  ++(-1,0) node [kropka] {}  ++(-1,0) node [kropka] {}   ++(-1,0) node [kropka] {}  ++(-1,0) node [kropka] {} ;

\path (-2,6) node [kropka] {}  ++(-1,0) node [kropka] {}    ++(-1,0) node [kropka] {}    ++(-1,0) node [kropka] {}   ++(-1,0) node [kropka] {} ;

\path (-2,7) node [kropka] {}   ++(-1,0) node [kropka] {}     ++(-1,0) node [kropka] {}     ++(-1,0) node [kropka] {}  ++(-1,0) node [kropka] {} ;

\path (-2,8) node [kropka] {}  ++(-1,0) node [kropka] {}    ++(-1,0) node [kropka] {}     ++(-1,0) node [kropka] {}   ++(-1,0) node [kropka] {} ;

\path (-2,9) node [kropka] {}   ++(-1,0) node [kropka] {}   ++(-1,0) node [kropka] {}     ++(-1,0) node [kropka] {}   ++(-1,0) node [kropka] {} ;

\path (-2,10) node [kropka] {}   ++(-1,0) node [kropka] {}   ++(-1,0) node [kropka] {}     ++(-1,0) node [kropka] {}   ++(-1,0) node [kropka] {} ;

\path (-2,11) node [kropka] {}   ++(-1,0) node [kropka] {}   ++(-1,0) node [kropka] {}     ++(-1,0) node [kropka] {}   ++(-1,0) node [kropka] {} ;

\path (-2,12) node [kropka] {}   ++(-1,0) node [kropka] {}   ++(-1,0) node [kropka] {}     ++(-1,0) node [kropka] {}   ++(-1,0) node [kropka] {} ;

\path (3,-8) node [kropka]{} ;

\path (3,-7) node [kropka]{} ;

\path (3,-6) node [kropka]{}  ++(-1,0)node [kropka]{} ;

\path (3,-5)  node [kropka]{}  ++(-1,0)node [kropka]{}  ++(-1,0)node [kropka]{}  ++(-1,0)node [kropka]{} ;

\path (3,-4) node [kropka]{}  ++(-1,0) node [kropka]{}  ++(-1,0) node [kropka]{}  ++(-1,0)node [kropka]{}  ++(-1,0) node [kropka] {} ;

\path (3,-3) node [kropka]{}  ++(-1,0) node [kropka]{}  ++(-1,0) node [kropka]{}  ++(-1,0)  node [kropka]{}  ++(-1,0) node [kropka] {} ;

\path (3,-2)  node [kropka]{} ++(-1,0)node [kropka]{} ++(-1,0)node [kropka]{}  ++(-1,0)node [kropka]{}  ++(-1,0) node [kropka] {}  ;

\path (3,-1) node [kropka]{}  ++(-1,0) node [kropka]{} ++(-1,0)node [kropka]{} ++(-1,0)node [kropka]{} ++(-1,0) node [kropka] {} ;

\path (3,0)  node [kropka]{}  ++(-1,0)node [kropka]{}  ++(-1,0)node [kropka]{}  ++(-1,0)node [kropka]{} ++(-1,0) node [kropka] {}   ;

\path (3,1) node [kropka]{} ++(-1,0) node [kropka]{}  ++(-1,0)node [kropka]{}  ++(-1,0)node [kropka]{}  ++(-1,0) node [kropka] {}  ;

\path (3,2)  node [kropka]{}  ++(-1,0)node [kropka]{}  ++(-1,0)node [kropka]{}  ++(-1,0)node [kropka]{}  ++(-1,0) node [kropka] {} ;

\path (3,3) node [kropka]{}  ++(-1,0) node [kropka]{}  ++(-1,0)node [kropka]{}  ++(-1,0)node [kropka]{} ++(-1,0) node [kropka] {} ;

\path (3,4) node [kropka]{} ++(-1,0) node [kropka]{} ++(-1,0)node [kropka]{} ++(-1,0)node [kropka]{} ++(-1,0) node [kropka] {}  ;

\path (3,5) node [kropka]{} ++(-1,0) node [kropka]{}  ++(-1,0) node [kropka]{}  ++(-1,0)  node [kropka]{} ++(-1,0)node [kropka] {}  ;

\path (3,6) node [kropka]{} ++(-1,0)node [kropka]{} ++(-1,0)node [kropka]{} ++(-1,0)node [kropka]{} ++(-1,0) node [kropka] {}    ;

\path (3,7) node [kropka]{} ++(-1,0)node [kropka]{}++(-1,0)node [kropka]{} ++(-1,0)node [kropka]{}  ++(-1,0) node [kropka] {}  ;

\path (3,8) node [kropka]{} ++(-1,0)node [kropka]{} ++(-1,0)node [kropka]{} ++(-1,0)node [kropka]{}  ++(-1,0) node [kropka] {}  ;

\path (3,9) node [kropka]{} ++(-1,0)node [kropka]{} ++(-1,0)node [kropka]{} ++(-1,0)node [kropka]{} ++(-1,0) node [kropka] {} ;

\path (3,10) node [kropka]{} ++(-1,0)node [kropka]{} ++(-1,0)node [kropka]{} ++(-1,0)node [kropka]{} ++(-1,0) node [kropka] {} ;

\path (3,11) node [kropka]{} ++(-1,0)node [kropka]{} ++(-1,0)node [kropka]{} ++(-1,0)node [kropka]{} ++(-1,0) node [kropka] {} ;

\path (3,12) node [kropka]{} ++(-1,0)node [kropka]{} ++(-1,0)node [kropka]{} ++(-1,0)node [kropka]{} ++(-1,0) node [kropka] {} ;

\path (-7,4) node [kropka] {} ;

\path (-7,5) node [kropka] {} ;

\path (-7,6) node [kropka] {} ++(-1,0) node [kropka] {} ;

\path (-7,7) node [kropka] {}   ++(-1,0) node [kropka] {}  ++(-1,0) node [kropka] {}   ++(-1,0) node [kropka] {} ;

\path (-7,8) node [kropka] {}   ++(-1,0) node [kropka] {}   ++(-1,0) node [kropka] {}   ++(-1,0) node [kropka] {} ++(-1,0) node [kropka] {} ;

\path (-7,9) node [kropka] {} ++(-1,0) node [kropka] {}  ++(-1,0) node [kropka] {}  ++(-1,0) node [kropka] {} ++(-1,0) node [kropka] {} ;

\path (-7,10) node [kropka] {} ++(-1,0) node [kropka] {}  ++(-1,0) node [kropka] {}  ++(-1,0) node [kropka] {} ++(-1,0) node [kropka] {} ;

\path (-7,11) node [kropka] {} ++(-1,0) node [kropka] {}  ++(-1,0) node [kropka] {}  ++(-1,0) node [kropka] {} ++(-1,0) node [kropka] {} ;

\path (-7,12) node [kropka] {} ++(-1,0) node [kropka] {}  ++(-1,0) node [kropka] {}  ++(-1,0) node [kropka] {} ++(-1,0) node [kropka] {} ;

\path (-13,12) node [kropka] {} ++(1,0) node [kropka] {} ++(0,-1) node [kropka] {} ++(0,-1) node [kropka] {};

\end{tikzpicture}
\end{center}

 The definition of a pair of generalized powers was originally given in \cite[Definition 7.3]{BKS}, but in \cite[Definition 4.7]{BKPS} it uses a more convenient notation. For reader convenience we give a definition with notation similar to Remark \ref{diagram_equiv}.
\begin{definition}\label{gen_pow}
Let there be given a periodic diagram $J$ with period $J_0$ and integers $(m,n)$ and a unitary $\mathcal{U}\in B(\mathcal{H}).$
Define
$H:=\{\mathbf{v}=\{v_{(i,j)}\}_{(i,j)\in J_0}: v_{i,j}\in\mathcal{H}\}$
and
\begin{align*}V_{1}\mathbf{v}=\mathbf{x}\text{ where }&x_{(i,j)}=\left\{
                                              \begin{array}{ll}
                                                 0, & \hbox{if $i=0$ and  $(m-1,j-n)\notin J_0,$}\\
                                                \mathcal{U}v_{(m-1,j-n)}, & \hbox{if $i=0$ and $(m-1,j-n)\in J_0,$}\\
                                                v_{(i-1,j)}, & \hbox{if $i=1,\dots,m-1$.} \\
                                              \end{array}
                                            \right.\\
V_{2}\mathbf{v}=\mathbf{y}\text{ where }& y_{(i,j)}=\left\{
                                              \begin{array}{ll}
                                                v_{(i,j-1)}, & \hbox{if $(i,j-1)\in J_0$,} \\
                                                0, & \hbox{if $(i,j-1)\notin J$.}
                                              \end{array}\right.\end{align*}
We call $(V_1, V_2)$ a pair of generalized powers given by a periodic  diagram $J$ and a unitary $\mathcal{U}$.
\end{definition}

By \cite[Theorem 7.2]{BKS} generalized powers $(V_1, V_2)$ are compatible and $V_2^{*n}V_1^m$ is unitary - equal to the extension of $\mathcal{U}$ onto $H$ given by $V_2^{*n}V_1^m\mathbf{v}=\{\mathcal{U}v_{(i,j)}\}_{(i,j)\in J_0}$, for $\mathbf{v}=\{v_{(i,j)}\}_{(i,j)\in J_0}$. 
Conversely, if $(V_1, V_2)$ is a commuting, compatible pair of isometries, such that $V_2^{*n}V_1^m$ is unitary for some $(m,n)$ then the pair may be decomposed among pairs of unitary operators and pairs of generalized powers.

Let us emphasis that the subspace $H_d$ in Theorem \ref{dec_thm} is not the maximal one, where a pair may be decomposed among pairs given by diagrams. To be precise, some special cases of pairs in other parts may be also considered as pairs given by a diagram. The simplest example is $M_{z_1}, M_{z_2}$ on $L^2(\mathbb{T})$ which is a pair of unitary operators but also a pair given by the diagram $\mathbb{Z}^2$. It turns out that some generalized powers may be also described as pairs given by diagrams. It is explained in \cite[p.176, Remark 4.9]{BKPS} - a pair of generalized powers given by a (periodic) diagram $J$ and a bilateral shift $\mathcal{U}$ with the wandering subspace $\mathcal{H}_w$ may be also described as a pair of isometries given by the same diagram $J$ and the space $\mathcal{H}_w$. The reverse also holds. That is, if $V_1, V_2$ is a pair given by a periodic diagram and the space $\mathcal{H}$ then $V_2^{*n}V_1^m$ is a bilateral shift, for which $\mathcal{H}$ is wandering where $m,n$ are integers corresponding to the periodicity of the diagram. Such the pair may be also described as a pair of generalized powers given by the same diagram and a unitary $\mathcal{U}=V_2^{*n}V_1^m|_{\bigoplus_{i\in\mathbb{Z}}(V_2^{*n}V_1^m)^i\mathcal{H}}$. 

The evanescent part of a random field with irrational half-plane past  is described by a continuously given pair by Theorem \ref{Thmirs}. The original definition of a continuously given pair \cite[Definition 1]{BKPSCorr} uses the concept of the shape of the diagram. To avoid recalling some technical construction we have formulated equivalent definition. 
\begin{definition}
A pair of compatible isometries $(V_1, V_2)\in B(H)$ is continuously given if the minimal $(V_1, V_2)$ reducing subspace containing $\ker V_1^*\cap\ker V_2^*$ is the whole $H$ and there is no nontrivial subspace reducing $(V_1, V_2)$ to a pair of generalized powers or a pair given by a diagram.
\end{definition}

\section{The evanescent part for the past given by $L$}\label{SectL}
In this section we consider the random field $\{X_{(s,t)}\}_{(s,t)\in \Zzz^2}$ with no deterministic part and the past given by the half-plane $L$ \eqref{QL}. Our aim is to decompose such a random field between a purely non-deterministic field and an evanescent field and show that the latter one consists of a one-dimensional deterministic and a one-dimensional completely non-deterministic processes. Such a result is not new, however we use operator theory approach, which applied in the next section will describe evanescent part of a random field with a rational half-plane past in the new way.

Let $$-L=\{(k,l): (-k,-l)\in L\}=\mathbb{Z}\times\mathbb{Z}_+\cup \mathbb{Z}_+\times\{0\}$$ and $V^{(k,l)}$ be as in Definition \ref{defSuciu} which in the case of $(k,l)\in -L,$ so $l\ge 0$ simplifies to $$V^{(k,l)}=\left\{\begin{array}{l}
                                                                                                   V_1^{k}V_2^{l} \text{ for } k\geq 0, \\
                                                                                                   V_1^{*-k}V_2^{l} \text{ for } k<0.
                                                                                               \end{array}
                                                                                                \right.$$

\begin{theorem}\label{Thm_HtHe}
Let $\{X_{(s,t)}\}_{(s,t)\in \Zzz^2}$  be a random field with no deterministic part and the half-plane past $L$ \eqref{QL}.
Further, let $U_1, U_2$ be unitaries corresponding to the random field as in \eqref{processToOp} and $V_1=U_1|_{H^{(0,0)}}, V_2=U_2|_{H^{(0,0)}}$ be isometries defined on $H^{(0,0)}$ - the past of $X_{(0,0)}$.

There is a decomposition $H^{(0,0)}=H_t^{(0,0)}\oplus H_e^{(0,0)}$ such that:

\begin{itemize}
\item $H_t^{(0,0)}, H_e^{(0,0)}$ reduce $V_1, V_2$,
\item $V_1|_{H_t^{(0,0)}}, V_2|_{{H_t}^{(0,0)}}$ is a totally non unitary pair,
\item $V_1|_{H_e^{(0,0)}}$ is unitary, $V_2|_{H_e^{(0,0)}}$ is a unilateral shift.
\end{itemize}

Moreover,
\begin{itemize}
\item there is a vector $v\in\ker V_2^*\cap H_e$ cyclic for $V_1|_{\ker V_2^*\cap H_e}$ and $$H_e=\bigoplus_{n\ge 0}\spa\{V_2^nV_1^mv: m=0,\dots,K-1\}$$ where $K=\dim\ker V_2^*\leq\infty,$
\item $V_1|_{H_t^{(0,0)}}, V_2|_{H_t^{(0,0)}}$ is given by the diagram $-L\cup\{(0,0)\}$ and $$H_t^{(0,0)}=\bigoplus_{(k,l)\in -L\cup\{(0,0)\}} V^{(k,l)}\ker V_1^*.$$
\end{itemize}
\end{theorem}
\begin{proof}
Since the random field $\{X_{(s,t)}\}_{(s,t)\in\Zzz^2}$ has no deterministic part, $V_1, V_2$ is completely non-unitary, so $V_1V_2$ is a unilateral shift, and we have $H^{(0,0)}=\bigoplus_{n\ge 0}(V_1V_2)^n\mathcal{E}$ where $\mathcal{E}=\ker (V_1V_2)^*.$ For the matrix representation of $V_1, V_2$ let it be given  $l^2(\mathcal{E})=\{(h_0,h_1,\dots ): h_n\in \mathcal{E}, \sum_{n\ge 0} \|h_i\|^2<\infty\}$ and a unitary operator
$$\Psi: H^{(0,0)} \ni h = \sum_{n=0}^\infty (V_1V_2)^nh_n \mapsto (h_0,h_1,\dots ) \in l^2(\mathcal{E}).$$ We show the decomposition $l^2(\mathcal{E})=H_t^W\oplus H_e^W$ satisfying conditions analogous to those in the statement for the pair $W_1=\Psi V_1 \Psi^*, W_2=\Psi V_2 \Psi^*\in B(l^2(\mathcal{E}))$. Then $H_t^{(0,0)}=\Psi^*H_t^W, H_e^{(0,0)}=\Psi^*H_e^W$ satisfy the statement for $V_1, V_2$.

Berger, Coburn, Lebow model \cite{BCL} via \cite[(3-6)]{Pop} provides a matrix representation
$$W_1 =\begin{bmatrix}
U(I-P) & 0 & 0 & 0 &  \dots \\
UP & U(I-P) & 0 & 0 & \dots \\
0 & UP & U(I-P) & 0 & \dots \\
0 & 0 & UP & U(I-P) &  \dots \\
0 & 0 & 0 & UP &  \dots \\
\vdots & \vdots & \vdots &\vdots & \ddots
\end{bmatrix}$$
and
$$W_2 =\begin{bmatrix}
PU^* & 0 & 0 & 0  &\dots \\
(I-P)U^* & PU^* & 0  & 0 & \dots \\
0 & (I-P)U^* & PU^*  & 0 & \dots \\
0 & 0 & (I-P)U^* & PU^*  & \dots \\
0 & 0 & 0 & (I-P)U^*   & \dots \\
\vdots & \vdots &\vdots & \vdots &  \ddots
\end{bmatrix}$$
where $U, P\in B(\mathcal{E})$
\begin{equation}\label{Udef}U=(W_1(I-W_2W_2^*)+(I-W_1W_1^*)W_2^*)|_\mathcal{E}=(W_1(I-W_2W_2^*)+W_2^*)|_\mathcal{E}\end{equation}
is unitary and
$$P= W_2(I-W_1W_1^*)W_2^*|_\mathcal{E}$$ is an orthogonal projection.
The middle term in \eqref{Udef} is more convenient to calculate $U^*$  as follows \begin{equation}\label{ustar}U^*=P_\mathcal{E}((I-W_2W_2^*)W_1^*+W_2(I-W_1W_1^*))|_\mathcal{E}= (W_1^*+W_2(I-W_1W_1^*))|_\mathcal{E}.\end{equation}

Let us define 
$$H^W_t:=\bigoplus_{(k,l)\in L\cup\{(0,0)\}} W^{(k,l)}\ker W_1^*.$$
To see that $H^W_t$ is well defined we have to show that $W^{(k,l)}\ker W_1^*$ are pairwise orthogonal.
We do it by proving that:
\begin{enumerate}
    \item[(a)] $\ker W_1^*=EUP\mathcal{E},$ where $E$ is a natural embedding of $\mathcal{E}$ into $l^2(\mathcal{E})$,
    \item[(b)] $UP\mathcal{E}$ is wandering for $U$,
    \item[(c)] there is a one-to-one correspondence between $\{W^{(k,l)}\ker W_1^*: (k,l)\in-L\cup\{(0,0)\} \}$ and $\{(W_1W_2)^iEU^nP\mathcal{E} : i\ge 0, n\in\mathbb{Z}\}$.
\end{enumerate}

To see $(a)$ let us calculate that \begin{equation}\label{w1w1*t}W_1W_1^*=\begin{bmatrix}
I-UPU^* & 0 &  0 & \dots \\
0 & I & 0 & \dots \\
0 & 0 & I &  \dots \\
\vdots & \vdots &  \vdots & \ddots
\end{bmatrix},\quad
W_2W_2^* =\begin{bmatrix}
P & 0 & 0 &   \dots \\
0 & I & 0 &   \dots \\
0 & 0 & I &   \dots \\

\vdots & \vdots & \vdots  & \ddots
\end{bmatrix}\end{equation}
and so $P_{\ker W_1^*}=I-W_1W_1^*=\begin{bmatrix}
UPU^* & 0 & \dots \\
0 & 0 &  \dots \\
\vdots & \vdots &\ddots \
\end{bmatrix}=EUPU^*P_{\mathcal{E}}.$   Hence $$\ker W_1^*=EUPU^*P_{\mathcal{E}}l^2(\mathcal{E})=EUPU^*\mathcal{E}=EUP\mathcal{E}.$$

Next let us see $(b)$.
The total order given by the past $L$ implies
$W_2(l^2(\mathcal{E}))=\Psi V_2(H^{(0,0)})=\Psi H^{(0,-1)}\subset \Psi H^{(-n,0)}=\Psi V_1^n(H^{(0,0)})=W_1^n(l^2(\mathcal{E}))$ and in turn \begin{equation}\label{comp1t} W_2W_2^*\leq W_1^{n}W_1^{*n}\text{ for }n\ge 1.\end{equation}
Now, let us prove that $PU^nP=0$ for $n\ne 0$.

The proof is done by induction, it requires to prove simultaneously formulae on $W_1^nW_1^{*n}$. 
More precisely, we show that 
$$W_1^nW_1^{*n}=\begin{bmatrix}
I-\sum_{i=1}^n U^iPU^{*i} & 0 & 0 &  \dots \\
0 & I & 0 &   \dots \\
0 & 0 & I &   \dots \\
\vdots & \vdots &\vdots & \ddots \
\end{bmatrix} \quad \textnormal{ and } \quad PU^nP=0.$$

Base step, $n=1$ is the result of \eqref{w1w1*t} and \eqref{comp1t} by which
$$P(I-UPU^*)=P \Rightarrow PUPU^*=0 \Rightarrow  PUP=0.$$

Inductive step: assume that $$W_1^nW_1^{*n}=\begin{bmatrix}
I-\sum_{i=1}^n U^iPU^{*i} & 0 & 0 & \dots \\
0 & I & 0  & \dots \\
0 & 0 & I  & \dots \\
\vdots & \vdots & \vdots & \ddots \
\end{bmatrix} \quad \textnormal{ and } \quad PU^iP=0 \text{ for
 } i=1,\dots n.$$
 Taking advantage of  $PU^iP=0$ and $PU^{*i}P=0$ for $i=1,\dots n$ one can directly calculate  $$W_1^{n+1}W_1^{*n+1}=W_1(W_1^nW_1^{*n})W_1^*=\begin{bmatrix}
I-\sum_{i=1}^{n+1} U^iPU^{*i} & 0 & 0 & 0 & \dots \\
0 & I & 0 & 0  & \dots \\
0 & 0 & I & 0  & \dots \\
0 & 0 & 0 & I   & \dots \\
\vdots & \vdots & \vdots & \hdots
\end{bmatrix}.$$
Hence by \eqref{comp1t} we get
$$P(I-\sum_{i=1}^{n+1} U^iPU^{*i})=P\Rightarrow \sum_{i=1}^{n+1} PU^iPU^{*i}=0\Rightarrow PU^{n+1}PU^{*n+1}=0\Rightarrow PU^{n+1}P=0$$ where the last but one implication is by inductive assumption $PU^iP=0$ for $i=1,\dots,n$. 

Therefore, $PU^nP=0$, for $n>0$ (so for $n\neq 0$) yields that $UP\mathcal{E}$ is wandering for $U$.

Now, let us focus on $(c)$.
By a direct calculation on matrix form of $W_1, W_2$ taking advantage of $PU^nP=0$, for $n\not=0$ one can check that
\begin{equation}\label{UbyW}
EU^{n}|_{UP\mathcal{E}}=W_1^{n}|_{EUP\mathcal{E}}=W_1^n|_{\ker W_1^*},\quad 
\end{equation} and
\begin{equation}\label{UbyW2}
EU^{*n}|_{UP\mathcal{E}}=EU^{*n}UPU^*|_{UP\mathcal{E}}=EU^{*(n-1)}PU^*|_{UP\mathcal{E}}=W_1^{*n-1}W_2|_{EUP\mathcal{E}}=W_1^{*n-1}W_2|_{\ker W_1^*},\end{equation}
 for $n\geq 1$.
Hence it follows the one-to-one correspondence:

$$(W_1W_2)^{l}EU^{k-l}(UP\mathcal{E})=W_1^{l}W_2^{l}W_1^{k-l}(\ker W_1^*)=W_1^{k}W_2^{l}(\ker W_1^*)=W^{(k,l)}(\ker W_1^*)$$ for $k\geq l \geq 0$ and \begin{align*}(W_1W_2)^{l-1}EU^{*l-k}(UP\mathcal{E})&=W_2^{l-1}W_1^{l-1}W_1^{*l-k-1}W_2(\ker W_1^*)\\&=\left\{
                                                                                                \begin{array}{l}
                                                                                                W_2^{l-1}W_1^{l-1}W_1^{*l-1}W_1^{k}W_2(\ker W_1^*) \text{ for } k\geq 0,\\
                                                                                                W_2^{l-1}W_1^{l-1}W_1^{*l-1}W_1^{*-k}W_2(\ker W_1^*) \text{ for } k<0
                                                                                                \end{array}
                                                                                              \right.
\\&=\left\{
                                                                                                \begin{array}{l}
                                                                                                W_2^{l-1}W_1^{k}W_2(\ker W_1^*) \text{ for } k\geq 0,\\
                                                                                                W_2^{l-1}W_1^{*-k}W_2(\ker W_1^*) \text{ for } k<0
                                                                                                \end{array}
                                                                                              \right.
\\&=\left\{
                                                                                                \begin{array}{l}
                                                                                                W_1^{k}W_2^{l}(\ker W_1^*) \text{ for } k\geq 0,\\
                                                                                                W_1^{*-k}W_2^{l}(\ker W_1^*) \text{ for } k<0
                                                                                                \end{array}
                                                                                              \right.=W^{(k,l)}(\ker W_1^*)\end{align*} for $k<l, l>0$ where  $W_1^{l-1}W_1^{*l-l}$ can be reduced, as it is a projection and factor in the product of contractions factorizing the initial isometry.  For similar reason $W_2^{l-1}W_1^{*-k}W_2=W_1^{*-k}W_1^{-k}W_2^{l-1}W_1^{*-k}W_2=W_1^{*-k}W_2^{l-1}W_1^{-k}W_1^{*-k}W_2=W_1^{*-k}W_2^{l}$ holds on $\ker W_1^*$.

Using analogous arguments one can check that $$W_1W^{(k,l)}(\ker W_1^*)=W^{(k+1,l)}(\ker W_1^*),\; W_2W^{(k,l)}(\ker W_1^*)=W^{(k,l+1)}(\ker W_1^*),$$ $$W_1^*W^{(k,l)}(\ker W_1^*)=\left\{\begin{array}{l} W^{(k-1,l)}(\ker W_1^*) \text{ if }(k-1,l)\in-L,\\\{0\} \text{ otherwise,}\end{array}\right.$$ $$W_2^*W^{(k,l)}(\ker W_1^*)=\left\{\begin{array}{l}W^{(k,l-1)}(\ker W_1^*)\text{ if }(k,l-1)\in-L,\\\{0\}\text{ otherwise,}\end{array}\right.$$ by which $H^W_t=\bigoplus_{(k,l)\in -L\cup\{(0,0)\}} W^{(k,l)}\ker W_1^*$ reduces $W_1, W_2$ to a pair given by the diagram $-L$, so to a totally non-unitary pair.

Since $\ker W_1^*\subset H^W_t$ we get $H^W_e=l^2(\mathcal{E})\ominus H^W_t$ reduces $W_1$ to a unitary operator. Hence, since the product is a unilateral shift $W_2|_{H^W_e}$ has to be a unilateral shift.
Denote for convenience $W_{ie}=W_i|_{H^W_e}$ for $i=1,2.$ It is well known that since $W_{1e}$ is unitary it doubly commute with $W_{2e}$ and $\ker W_{2e}^*$ reduces $W_{1e}$. Let us calculate $\ker W_{2e}^*$.
Note that 
$$\ker W_{2e}^*=\spa\{(I-W_2W_2^*)P_{H_e^W}\Psi X_{(n,0)}: n< 0,(I-W_2W_2^*)P_{H_e^W}\Psi X_{(n,-1)}: n\geq0 \}.$$
However, 
\begin{align*}&P_{H_e^W}\Psi X_{(n,-1)}=P_{H_e^W}W_1^{*(n+1)}W_1^{n+1} \Psi X_{(n,-1)}=P_{H_e^W}W_1^{*(n+1)}\Psi V_1^{n+1} X_{(n,-1)}=\\&
W_1^{*(n+1)}P_{H_e^W} \Psi X_{(-1,-1)}=W_1^{*(n+1)}P_{H_e^W} \Psi V_2 X_{(-1,0)}=W_1^{*(n+1)}P_{H_e^W} W_2\Psi X_{(-1,0)}=\\&
W_1^{*(n+1)}W_2P_{H_e^W} \Psi X_{(-1,0)}=W_{1e}^{*(n+1)}W_{2e}P_{H_e^W} \Psi X_{(-1,0)}=\\&
W_{2e}W_{1e}^{*(n+1)}P_{H_e^W} \Psi X_{(-1,0)}=W_2W_1^{*(n+1)}P_{H_e^W} \Psi X_{(-1,0)}\end{align*}

and so $P_{H_e^W} \Psi X_{(n,-1)}\in W_2(H_e^W)$ for $n\geq 0$.
Hence $$\ker W_{2e}^*=\spa\{(I-W_2W_2^*)P_{H_e^W}\Psi X_{(n,0)}: n< 0\}.$$  Since $W_1(I-W_2W_2^*)P_{H_e^W}\Psi X_{(n,0)}=(I-W_2W_2^*)P_{H_e^W}\Psi X_{(n-1,0)}$ the vector $(I-W_2W_2^*)P_{H_e^W}\Psi X_{(-1,0)}$ is cyclic for $W_{1e}|_{\ker W^*_{2e}}$. The formulae on $H_e^W$ is immediate consequence.
\end{proof}

From the proof above we get the remark to be used in the next section.
\begin{remark}
Under assumptions of Theorem \ref{Thm_HtHe} and notations from its proof one can see that
\begin{equation}\label{wandering}W^{(i,j)}W^{(k,l)}(UP\mathcal{E})=W^{(k+i, l+j)}(UP\mathcal{E}), \textnormal{ for }  (i,j), (k,l)\in -L, \end{equation} and for $(k,l), (i,j)\in-L$ such that $(k-i, l-j)\in-L$
\begin{equation}\label{wandering2}W^{(i,j)*}W^{(k,l)}(UP\mathcal{E})=W^{(k-i, l-j)}(UP\mathcal{E}).\end{equation}
\end{remark}

Assume notation of Theorem \ref{Thm_HtHe}, however it will be more convenient to consider operators on $H^{(1,0)}$ where $X_{(0,0)}$ belongs instead on $H^{(0,0)}$, in other words $V_1=U_1|_{H^{(1,0)}}$ and $V_2=U_2|_{H^{(1,0)}}$. Since $H^{(0,0)}= U_1 H^{(1,0)}=V_1 H^{(1,0)}$ we get $\ker V_1^*\subset  H^{(1,0)}\ominus H^{(0,0)}$. Since $H^{(1,0)}=H^{(0,0)}\vee \{X_{(0,0)}\}$ we get $\ker V_1^*=\mathbb{C}e_{(0,0)}$ where $e_{(0,0)}=\frac{P_{\ker V_1^*} X_{(0,0)}}{\|P_{\ker V_1^*} X_{(0,0)}\|}$. Hence $\ker V_1^*$ is at most one-dimensional and $H_t^{(1,0)}\ne \{0\}$ if and only if $X_{(0,0)}\notin H^{(0,0)}$. Moreover, $e_{(k,l)}=V^{(k,l)}e_{(0,0)}$, for $(k,l)\in -L\cup\{(0,0)\}$, form orthonormal basis of $H_t^{(1,0)}$.

Let us denote for short $P_{H_t^{(1,0)}}=P_t^{(1,0)}$. Then $P_t^{(1,0)}  X_{(0,0)}=
\sum_{(k,l)\in -L\cup\{(0,0)\}} \alpha_{(k,l)} e_{(k,l)}.$
Let us consider $(-s,-t)\in L$ and calculate $I_{(-s,-t)}$ an innovation part of the vector $P_t^{(1,0)}X_{(-s,-t)}$. We get
$$P_t^{(1,0)}X_{(-s,-t)} = P_t^{(1,0)}U_1^sU_2^tX_{(0,0)} =P_t^{(1,0)}V_1^sV_2^t X_{(0,0)}= V_1^sV_2^t P_t^{(1,0)} X_{(0,0)}= \sum_{(l,k)\in -L\cup\{(0,0)\}} \alpha_{(k,l)} e_{(k+s,l+t)}.$$

Hence $ P_t^{(1,0)} H^{(-s,-t)}=\spa\{e_{(k+s,l+t)}:(k,l)\in-L\}$ and $$I_{(-s,-t)}= P_t^{(1,0)} X_{(-s,-t)}-P_{H^{(-s,-t)}}P_t^{(1,0)} X_{(-s,-t)}=\alpha_{(s,t)}e_{(s,t)}.$$ In other words $H_t^{(1,0)}$ is a subspace which reduces the random field $X_{(i,j)}$ to purely non-deterministic random field, where the innovation part $I_{(i,j)}$ is equal $\alpha_{(i,j)}e_{(-i,-j)},$ for $(i,j)\in L\cup \{(0,0)\}$.

\begin{theorem}\label{modelL}
The evanescent part $E_{(m,n)}$ of a weak-stationary random field with a half-plane past $L$ is deterministic along the first coordinate of the index and stochastic along the second coordinate that is, there is a random field $\{F_{(s,t)}\}_{(s,t)\in\mathbb{Z}^2}$ such that:
 $$ E_{(m,n)}=\sum_{t=0}^\infty \sum_{s=0}^{K-1}\beta_{(s,t)}F_{s+m,t+n},$$
and $\{F_{s_t,t}\}_{t\in\mathbb{Z}}$ is a one-dimensional white noise for any sequence $\{s_t\}$, and $\{F_{s,t}\}_{s=0}^{K-1}$ is deterministic process for any $t$, where $K=\dim\{F_{s,0}:s=0,1,\dots\}\leq\infty.$
\end{theorem}
\begin{proof}
The statement follows directly from Theorem \ref{Thm_HtHe}. Indeed, let $F_{(0,0)}$ be the cyclic vector and $F_{(s,t)}=U_1^{-s}U_2^{-t}F_{(0,0)}$ for $(s,t)\in\mathbb{Z}^2$.
Properties of $\{F_{(s,t)}\}_{(s,t)\in\mathbb{Z}^2}$ follows directly from types of $V_1, V_2$.
\end{proof}

\begin{example}
  Let $\{\alpha_n\}_{n=-\infty}^\infty$ be a sequence of orthogonal (i.e. $\mathbb{E}\alpha_n\overline{\alpha_m}=0$ for $m\not=n$) random variables.
  Let $X_{(s,t)}:=\alpha_t$. Then $\{X_{(s,t)}\}_{s,t\in\Zzz}$ is a weak-stationary random field.

For the past given by $L$ we have $H^{(s,t)}=\bigvee\{X_{(i,j)}:(i,j)\in L + \{(s,t)\}\}=\bigvee\{\alpha_i: i\leq t\}$. Thus $X_{(s,t)}\in H^{(s,t)}$ and $\bigcap\limits_{(s,t)\in \Zzz^2} H^{(s,t)} = \{0\}$. Hence the random  field $\{X_{(s,t)}\}_{s,t\in\Zzz}$ is evanescent. Moreover, $H^{(s,t)}_X$
reduces $(U_1,U_2)$ to the pair: an identity, a unilateral shift of multiplicity $1$.
\end{example}

\section{The rational nonsymmetrical half-plane}\label{SectRS}

The half-plane $S_{\textbf{v}}$ defined in \eqref{hp}, in the case $\textbf{v}\in \Zzz^2$ is called the rational nonsymmetrical half-plane (see e.g. \cite{Cuny,FMP}). In other words, rational nonsymmetrical half-plane $S_{\textbf{v}}$ is a half-plane with a rational slope.  Note that $L$ \eqref{QL} is a half-plane with a slope $0$, so the special case of the rational slope. It turns out that an arbitrary rational nonsymmetrical half-plane is a rotation of $L$ and via such rotation the result of previous section for the past given by $L$ can be applied to a random field with the past given by an arbitrary rational nonsymmetrical half-plane. 
\begin{remark}\label{psi_rem}
For $\textbf{v}=(k,l)\in\Zzz^2_-$ where $k,l$  are relatively prime, there are $(p,q)\in \Zzz^2$ such that $pk+ql=-1$ and $-l>p>0\geq q>k$ and
$$\psi:L\ni(m,n)\mapsto \left[\begin{array}{cc}-l &p\\k &q\end{array}\right]\left[\begin{array}{c}m\\n\end{array}\right]= (pn-lm,qn+km)\in S_{(k,l)}$$ is a bijection.\end{remark}
\begin{proof}Since $-k, -l$ are positive, relatively prime, remainders of $\frac{-k}{-l}, \frac{-2k}{-l},\dots,\frac{-(-l-1)k}{-l}$ are positive and pairwise different. Hence the set of such remainders is the whole $\{1,\dots, -l-1\},$ so there is $p\in\{1,\dots,-l-1\}$ such that the remainder of $\frac{-pk}{-l}$ is $1$. Let $q'$ be the corresponding quotient, so $0\leq q'< -k$. Consequently $-pk=-q'l+1$ and in turn $pk-q'l=-1$. Taking $q=-q'$ we get $pk+ql=-1$ and $-l>p>0\geq q>k$.
Hence $\psi$ may be defined on the whole $\mathbb{Z}^2$, where it is linear and invertible. Consequently, it is enough to check that $\psi(L)\subset S_{(k,l)}$ and $\psi^{-1}(S_{(k,l)})\subset L$. 

For $n=0, m<0$ we get $\seq{(-lm,km),(k,l)}=0$, while for $n<0, m\in\mathbb{Z}$
$\seq{(pn-lm,qn+km),(k,l)}= (pk+ql)n=-n>0$. Hence $\psi(L)\subset S_{(k,l)}$. Note that $$\psi^{-1}(a,b)=\left[\begin{array}{cc}-l &p\\k &q\end{array}\right]^{-1}\left[\begin{array}{c}a\\b\end{array}\right]=(qa-pb,-ka-lb).$$ For $(a,b)\in S_{(k,l)}$ either $-ka-lb<0$ and $a$ is arbitrary or $-ka-lb=0$ and $a<0$. For the first case clearly $(qa-pb,-ka-lb)\in L$. For the second case, since $b=\frac{-ka}{l}$ we get $qa-pb=qa-\frac{-pka}{l}=\frac{a}{l}(ql+pk)=-\frac{a}{l}$ which by $a<0$ is less then $0$. Hence $(qa-pb,-ka-lb)\in\mathbb{Z}_-\times\{0\}\subset L$, so $\psi^{-1}(S_{(k,l)})\subset L$.\end{proof}

As we mentioned earlier, version of Theorem \ref{Thm_HtHe} for any rational nonsymmetrical half-plane follows from Remark \ref{psi_rem}.

\begin{theorem}\label{Thmrs}
Let $\{X_{(s,t)}\}_{(s,t)\in \Zzz^2}$  be a random field with no deterministic part and the past given by the rational nonsymmetrical half-plane past $S_{\textbf{v}}$ \eqref{hp} for $\textbf{v}=(k,l)\in\mathbb{Z}^2$.
Further, let $U_1, U_2$ be unitaries corresponding to the random field as in \eqref{processToOp} and $V_1=U_1|_{H^{(0,0)}}, V_2=U_2|_{H^{(0,0)}}$ be isometries defined on $H^{(0,0)}$ - the past of $X_{(0,0)}$.

There is a decomposition $H^{(0,0)}=H_t^{(0,0)}\oplus H_e^{(0,0)}$ such that:

\begin{itemize}
\item $H_t^{(0,0)}, H_e^{(0,0)}$ reduce $V_1, V_2$,
\item $V_1|_{H_t^{(0,0)}}, V_2|_{{H_t}^{(0,0)}}$ is a totally non unitary pair,
\item $V_1|_{H_e^{(0,0)}}, V_2|_{H_e^{(0,0)}}$ is a pair of generalized powers given by a cyclic unitary.
\end{itemize}

Moreover, if $k,l\in\mathbb{Z}_-$ and are relatively prime  we get
\begin{itemize}
\item  $(V_1|_{H_e^{(0,0)}}, V_2|_{H_e^{(0,0)}})$ is a pair of generalized powers given by a diagram with period $J_0=\{(i,j)\in-S_{(k,l)}:0\leq i\leq -l-1\}\cup\{(0,0)\},$ integers $(-l,-k)$, the space $\mathcal{H}=\ker V_1^{*p}V_2^{-q}\cap H_e^{(0,0)}$ where $p,q$ are as in Remark \ref{psi_rem} and a unitary $\mathcal{U}=V_2^{*-k}V_1^{-l}|_{\ker V_1^{*p}V_2^{-q}}$, and $$H_e^{(0,0)}=\bigoplus_{(i,j)\in J_0} \bigvee\{(V_2^{*-k}V_1^{-l})^\kappa V^{(i,j)}w:\kappa\ge 0\},$$ where $w$ is the cyclic vector for $\mathcal{U}$ (see Definition \ref{gen_pow}),
\item $(V_1|_{H_t^{(0,0)}}, V_2|_{H_t^{(0,0)}})$ is given by the diagram $-S_{\textbf{v}}\cup\{(0,0)\}$ and $$H_t^{(0,0)}=\bigoplus_{(i,j)\in -S_{(k,l)}\cup\{(0,0)\}} V^{(i,j)}\ker V_1^{*-l}V_2^{-k}.$$
\end{itemize}
\end{theorem}
\begin{proof}
Without loss of generality we may assume since the beginning that $\textbf{v}=(k,l)\in\Zzz^2_-$ where $k,l$  are relatively prime as in Remark \ref{psi_rem} and use $\psi$ defined there. By linearity of $\psi$ there is a natural extension of $\psi$ on $\mathbb{Z}^2$ given by $\psi(i,j)=-\psi(-i,-j)$ for $(i,j)\in-L$ and $\psi(0,0)=(0,0)$. We use this extension if needed.
Let us define a random field $\{Y_{(m,n)}\}_{(m,n)\in \Zzz^2}$ by $$Y_{(m,n)}:=X_{\psi(m,n)}=X_{(pn-lm,qn+km)}.$$ Since formulae in $\psi$ are linear, the cross moments of $\{Y_{(m,n)}\}_{(m,n)\in \Zzz^2}$ inherit dependence only on the distance and so  $\{Y_{(m,n)}\}_{(m,n)\in \Zzz^2}$ is weak-stationary. Note that definition of $\psi$ implies that for any $(m,n)\in\mathbb{Z}^2$ the past spaces $H^{(m,n)}$ given by \eqref{Hij} are spanned by the same random variables when calculated for the random field $\{Y_{(m,n)}\}_{(m,n)\in \Zzz^2}$ with the past $L$ as well as for the random field $\{X_{(s,t)}\}_{(s,t)\in \Zzz^2}$ with the past $S_{\textbf{v}}$. In other words the past of $X_{(m,n)}$ and the past of $Y_{(m,n)}$ is the same space   $H^{(m,n)}$. On the other hand, operators $U^Y_1, U^Y_2$ defined by \eqref{processToOp} for the random field $\{Y_{(m,n)}\}_{(m,n)\in \Zzz^2}$ are different from $U_1, U_2$. However, there is a connection between $(U^Y_1, U^Y_2)$ and $(U_1, U_2)$:
$$U^Y_1Y_{(m,n)}=Y_{(m-1,n)}=X_{(pn-lm+l,qn+km-k)}=U_2^{*-k}U_1^{-l}X_{(pn-lm,qn+km)}=U_2^{*-k}U_1^{-l}Y_{(m,n)},$$
$$U^Y_2Y_{(m,n)}=Y_{(m,n-1)}=X_{(pn-lm-p,qn+km-q)}=U_2^{*-q}U_1^{p}X_{(pn-lm,qn+km)}=U_2^{*-q}U_1^{p}Y_{(m,n)}$$
and by similar calculation $U_1=U_1^{Y*-q}U_2^{Y-k}, U_2=U_1^{Y*p}U_2^{Y-l}$.
Consequently $(U_1, U_2)$ and $(U_1^Y, U_2^Y)$ have the same reducing subspaces and the same subspaces reduce them to unitary pairs. Thus $\{Y_{(i,j)}\}_{(s,t)\in \Zzz^2}$ has no deterministic part, as it was assumed for $\{X_{(i,j)}\}_{(s,t)\in \Zzz^2}$ and we may apply Theorem \ref{Thm_HtHe} to $V_1^Y=U_1^Y|_{H^{(0,0)}}, V_2^Y=U_2^Y|_{H^{(0,0)}}.$

Let $H^{(0,0)}=H_t^{Y(0,0)}\oplus H_e^{Y(0,0)}$ be the decomposition given by Theorem \ref{Thm_HtHe} for $(V_1^Y, V_2^Y)$. Let us check that $H_t^{Y(0,0)}, H_e^{Y(0,0)}$ reduce $(V_1, V_2)$. Since the subspaces complements each other to $H^{(0,0)}$ and $V_1, V_2\in B(H^{(0,0)})$ it is enough to show that $H_t^{Y(0,0)}, H_e^{Y(0,0)}$ are $(V_1, V_2)$ invariant. Let $H_t, H_e$ be the minimal $(U_1^Y, U_2^Y)$ reducing subspaces containing $H_t^{Y(0,0)}, H_e^{Y(0,0)}$ respectively. Hence they are $(U_1, U_2)$ reducing. On the other hand, by Remark \ref{projprop} $H_t^{Y(0,0)}=H_t\cap H^{(0,0)}, H_e^{Y(0,0)}=H_e\cap H^{(0,0)}$ which clearly are $U_1, U_2$ invariant subspaces of $H^{(0,0)}$, so they are $V_1, V_2$ invariant.

Next we show that
$H_t^{Y(0,0)}=H_t^{(0,0)}$ and so $H_e^{Y(0,0)}=H_e^{(0,0)}$ and that restrictions of $(V_1, V_2)$ are as described in the statement.

Let us first  show that $H_t^{Y(0,0)}$ and $H_t^{(0,0)}$ are sums of the same subspaces, but differently indexed and so they are equal. Note that $V_1^Y=U_1^Y|_{H^{(0,0)}}=U_2^{*-k}U_1^{-l}|_{H^{(0,0)}}=P_{H^{(0,0)}}U_2^{*-k}U_1^{-l}|_{H^{(0,0)}}=V_2^{*-k}V_1^{-l}$ and similarly $V_2^Y=V_2^{*-q}V_1^{p}.$ Hence $\ker {V^Y_1}^*=\ker V_1^{*-l}V_2^{-k}$ and $(V^Y)^{(i,j)}=V^{\psi(i,j)}$. 
\begin{align*}H_t^{Y(0,0)}=\bigoplus_{(i,j)\in -L\cup\{(0,0)\}} V^{Y(i,j)}\ker V_1^{Y*}=\bigoplus_{(i,j)\in -L\cup\{(0,0)\}} V^{\psi(i,j)}\ker V_1^{*-l}V_2^{-k}\\=\bigoplus_{(i,j)\in -S_{\textbf{v}}\cup\{(0,0)\}}V^{(i,j)}\ker V_1^{*-l}V_2^{-k}=H_t^{(0,0)}.\end{align*} Hence we may identify $V_1|_{H_t^{(0,0)}}, V_2|_{H_t^{(0,0)}}$ as a pair given by a diagram $-S_{\textbf{v}}\cup\{(0,0)\}$ and $\mathcal{H}=\ker V_1^{*-l}V_2^{-k}$ as in Remark \ref{diagram_equiv}.

Let us check that the pair $(V_1|_{H_e^{(0,0)}},V_2|_{H_e^{(0,0)}})$ is a pair of generalized powers. By Theorem \ref{Thm_HtHe} we have $V_1^{Y}|_{H_e^{(0,0)}}=V_2^{*-k}V_1^{-l}|_{H_e^{(0,0)}}$ a unitary, and $v$ a cyclic vector for $V_2^{*-k}V_1^{-l}|_{\ker V_2^{Y*}\cap{H_e^{(0,0)}}}$. Denote $\mathcal{H}_{(i,j)}=\bigvee\{(V_2^{*-k}V_1^{-l})^\kappa V^{(i,j)}v:\kappa>0\}$ for $(i,j)\in J_0$, so in particular $\mathcal{H}_{(0,0)}=\ker V_2^{Y*}\cap{H_e^{(0,0)}}$ and it reduces $V_2^{*-k}V_1^{-l}$. Moreover, $V_2^{*-k}V_1^{-l}|_{H_e^{(0,0)}}$ as a unitary doubly commutes with $V_1|_{H_e^{(0,0)}}, V_2|_{H_e^{(0,0)}}$ and in turn it commutes with each $V^{(i,j)}|_{H_e^{(0,0)}}$. Hence $\mathcal{H}_{(i,j)}=V^{(i,j)}(\ker V_2^{Y*}\cap{H_e^{(0,0)}})=V^{Y\psi^{-1}(i,j)}(\ker V_2^{Y*}\cap{H_e^{(0,0)}})=V_2^{Y\psi^{-1}(i,j)_2}V_1^{Y\psi^{-1}(i,j)_1}(\ker V_2^{Y*}\cap{H_e^{(0,0)}})=V_2^{Y\psi^{-1}(i,j)_2}(\ker V_2^{Y*}\cap{H_e^{(0,0)}})$ where $\psi^{-1}(i,j)=(\psi^{-1}(i,j)_1, \psi^{-1}(i,j)_2).$  One can check that $J_0\ni (i,j)\mapsto \psi^{-1}(i,j)_2\in\mathbb{Z}_+\cup\{0\}$ is injective (even bijective) and so $\mathcal{H}_{(i,j)}$ are pairwise orthogonal. Moreover, since each $V^{(i,j)}$ is an isometry for $(i,j)\in-S_{\textbf{v}}$ we get $V_1\mathcal{H}_{(i,j)}=\mathcal{H}_{(i+1,j)}$ for $i\leq-l-1$ and $V_1\mathcal{H}_{(-l-1,j)}=V_1^{-l}\mathcal{H}_{(0,j)}=V_1^{-l}V_2^{-k}V_1^{*-l}\mathcal{H}_{(0,j)}=\mathcal{H}_{(0,j-k)}$ and $V_2\mathcal{H}_{(i,j)}=\mathcal{H}_{(i,j+1)}$, so $\bigoplus_{(i,j)\in J_0} \mathcal{H}_{(i,j)}$ reduces $(V_1, V_2)$. Hence, in particular $\bigoplus_{(i,j)\in J_0} \mathcal{H}_{(i,j)}$ reduces a unilateral shift $V_2^Y$ and as it contains $\ker V_2^{Y*}\cap{H_e^{(0,0)}}$ we get $H_e^{(0,0)}=\bigoplus_{(i,j)\in J_0} \mathcal{H}_{(i,j)}.$  Consequently we may identify $V_1|_{H_e^{(0,0)}}, V_2|_{H_e^{(0,0)}}$ as a pair of generalized powers given by $J_0$, integers $(-l,-k)$, the space $\mathcal{H}=\ker V_2^{Y*}\cap{H_e^{(0,0)}}$  and the unitary $\mathcal{U}=V_2^{*-k}V_1^{-l}|_{\ker V_2^{Y*}\cap{H_e^{(0,0)}}}$ as in Definition \ref{gen_pow}.

\end{proof}

As a corollary we get a nice description of an evanescent part for the considered past.
\begin{theorem}\label{proprs}
    Let $\{E_{(s,t)}\}_{(s,t)\in \Zzz^2}$ be an evanescent part of a weak-stationary random field with a half-plane past $S_{\textbf{v}}$, where $\textbf{v}=(k,l)\in \Zzz_-^2$ for $k,l$ relatively prime. Then
there is a random field $\{F_{(s,t)}\}_{(s,t)\in \mathbb{Z}^2}$ such that:
 $$ E_{(m,n)}=\sum_{(s,t)\in J_{0,\dots,K-1}} \beta_{(s,t)}F_{s+m,t+n},$$
where $\{F_{m-li,n-ki}\}_{i\in\mathbb{Z}}$ is deterministic stochastic process for any $(m,n)$ and  $\{F_{m_i,n_i}\}_{i\in\mathbb{Z}}$ is a one-dimensional white noise whenever $(m_i-m_j,n_i-n_j)$ is not a multiplicity of $(-l,-k)$ for any $i\ne j$ where $J_{0,\dots,K-1}=\{(s,t)\in-S_{\textbf{v}}\cup\{(0,0)\}: 0\leq s\leq -lK-1\}$ and $K=\dim\{F_{-li,-ki}:i=0,1,\dots\}.$
\end{theorem}

\begin{proof}
The result is a consequence of Theorem \ref{Thmrs} on the space $H^{(1,0)}$. Indeed, take $F_{(0,0)}$ as the cyclic vector and define $F_{(s,t)}=U_1^sU_2^t F_{(0,0)}$. Moreover, we get $\{(V_2^{*-k}V_1^{-l})^\kappa V^{(i,j)}F_{(0,0)}:\kappa\ge 0\}=\{(V_2^{*-k}V_1^{-l})^\kappa V^{(i,j)}F_{(0,0)}:\kappa=0,\dots, K-1\}$. Hence $E_{(0,0)}\in H_e^{(1,0)}$ is a linear combination of $(V_2^{*-k}V_1^{-l})^\kappa V^{(i,j)}F_{(0,0)}=V^{(i-\kappa l,j+\kappa k)}F_{(0,0)}$ for $\kappa=0,\dots, K-1$ and $(i,j)\in J_0$. One can check that $\{(i-\kappa l,j+\kappa k):(i,j)\in J_0, \kappa=0,\dots, K-1\}=\bigcup_{\kappa=0}^{K-1}J_0+\kappa(-l,k)=J_{0,\dots,K-1}.$ Hence $E_{(0,0)}=\sum_{(s,t)\in J_{0,\dots,K-1}} \beta_{(s,t)}F_{s,t}$ for some coefficients $\beta_{(s,t)}$. Since $E_{(m,n)}=U_1^mU_2^nE_{(0,0)}$ we get the result by the definition of $F_{(s,t)}$.
\end{proof}

\begin{example}\label{exgpev}
Let $\{\alpha_n\}_{n=-\infty}^\infty$ be a sequence of orthogonal (i.e. $\mathbb{E}\alpha_n\overline{\alpha_m}=0$ for $m\not=n$) random variables.
Let $X_{(s,t)}:=\alpha_{s+t}$. Then $\{X_{(s,t)}\}_{s,t\in\Zzz}$ is a weak-stationary random field.
Let the past be  the half-plane $S_{[-1,-1]}$.
We have $H^{(s,t)}=\bigvee\{X_{(i,j)}:(i,j)\in S_{\bf{v}} + \{(s,t)\}\}=\bigvee\{\alpha_i: i\leq s+t \}$. Thus $X_{(s,t)}\in H^{(s,t)}$ and $\bigcap\limits_{(s,t)\in \Zzz^2} H^{(s,t)} = \{0\}$. Hence the random  field $\{X_{(s,t)}\}_{s,t\in\Zzz}$ is evanescent.

Moreover, we have $V_1|_{H^{(s,t)}} = V_2|_{H^{(s,t)}}$. So, in particular $(V_1|_{H^{(s,t)}},V_2|_{H^{(s,t)}})$ is a pair of generalized powers.
\end{example}
The following example shows the dependence of the type of random field on the past.
\begin{example}\label{exForS}
Let us consider the same weak-stationary random field as in Example \ref{exgpev} but with the past given by $L$.

Then we have $H^{(s,t)}=\bigvee\{X_{(i,j)}:(i,j)\in L + \{(s,t)\}\}=\bigvee\{\alpha_i: i\Zzz \}$ and so $\bigcap\limits_{(s,t)\in \Zzz^2} H^{(s,t)} = H$. Hence the random  field $\{X_{(s,t)}\}_{s,t\in\Zzz}$ is deterministic.
\end{example}

\section{The irrational nonsymmetrical half-plane}\label{SectNS}
In this section we consider the random field $\{X_{(s,t)}\}_{(s,t)\in \Zzz^2}$ with no deterministic part and the past given by the half-plane past $S_{\textbf{v}}$ where $v=(q,r)$ with $\frac{q}{r}\notin\mathbb{Q}$ that is the half-plane with irrational slope. Clearly, there can be constructed a bijection between $L$ and $S_{\textbf{v}}$  but not as a rotation as in Remark \ref{psi_rem}. Indeed, the rotation would have the angle corresponding to irrational slope, so it would not lead to $\mathbb{Z}^2$. In fact a bijection from $L$ to irrational half-plane would not  preserve dependence on distance and in turn a random field $\{Y_{(i,j)}\}_{(i,j)\in\mathbb{Z}}$ as defined in the proof of Theorem \ref{Thmrs} would not be weak-stationary. 

\begin{theorem}\label{Thmirs}
Let $\{X_{(s,t)}\}_{(s,t)\in \Zzz^2}$  be a random field with no deterministic part and the past given by $S_{\textbf{v}}$ \eqref{hp} where $v=(q,r)$ with $\frac{q}{r}\notin\mathbb{Q}$, that is the irrational nonsymmetrical half-plane past.
Further, let $U_1, U_2$ be unitary operators corresponding to the random field as in \eqref{processToOp} and $V_1=U_1|_{H^{(0,0)}}, V_2=U_2|_{H^{(0,0)}}$ be isometries defined on $H^{(0,0)}$ - the past of $X_{(0,0)}$.

There is a decomposition $H^{(0,0)}=H_t^{(0,0)}\oplus H_e^{(0,0)}$ such that:

\begin{itemize}
\item $H_t^{(0,0)}, H_e^{(0,0)}$ reduce $V_1, V_2$,
\item $V_1|_{H_t^{(0,0)}}, V_2|_{{H_t}^{(0,0)}}$ is a totally non unitary pair,
\item $V_1|_{H_e^{(0,0)}}, V_2|_{H_e^{(0,0)}}$ is a continuously given pair.
\end{itemize}
\end{theorem}
\begin{proof}
Since $\frac{q}{r}\notin\mathbb{Q}$ for any pair of integers $m,n$ we get $\left<(m,-n),\textbf{v}\right>\ne 0$.  Hence either $\left<(m,-n),\textbf{v}\right>< 0$ or $\left<(-m,n),\textbf{v}\right>< 0$  and one can show that \begin{equation}\label{irr}\{(i+\kappa m+k,j-\kappa n+l):(i,j)\in S_{\textbf{v}},\kappa\in\mathbb{Z}, k,l\ge 0\}=\mathbb{Z}^2.\end{equation}

Assume $\mathcal{H}\subset H^{(0,0)}$ is a subspace that reduces $(V_1, V_2)$ such that $V_2^{*n}V_1^m|_{\mathcal{H}}$ is unitary for some nonnegative integers $m,n$. Then $\mathcal{H}=\spa\{P_{\mathcal{H}}X_{(i,j)}:(i,j)\in S_{\textbf{v}}\}$ and by \eqref{irr} all $P_{\mathcal{H}}X_{(i,j)}$ have the same norm. Indeed, since $V_2^{*n}V_1^m|_{\mathcal{H}}$ is unitary we have $\|V_2^lV_1^k(V_2^{*n}V_1^m)^\kappa P_{\mathcal{H}}X_{(i,j)}\|=\|P_{\mathcal{H}}V_2^lV_1^k(V_2^{*n}V_1^m)^\kappa X_{(i,j)}\|=\|P_{\mathcal{H}}X_{(i+\kappa m+k,j-\kappa n+l)}\|$ and $P_{\mathcal{H}}X_{(i+\kappa m+k,j-\kappa n+l)}\in\mathcal{H}$ for $\kappa\in\mathbb{Z}$ and $k,l$ nonnegative integers. By \eqref{irr} we get $\mathcal{H}=\spa\{P_{\mathcal{H}}X_{(i,j)}:(i,j)\in\mathbb{Z}^2\}$. Hence $\mathcal{H}$ reduces $V_1, V_2$ to a unitary pair. Since the process is completely non deterministic we get $\mathcal{H}=\{0\}$. In other words we have shown that $V_2^{*n}V_1^m$ is not unitary  in restriction to any non-trivial reducing subspace.

Let us decompose $(V_1, V_2)$ by Theorem \ref{dec_thm}. Since $V_2^{*n}V_1^m$ is not unitary on subspaces reducing $(V_1, V_2)$ for $m,n$ non-negative integers, so also $m=0$ or $n=0$ we get $H_{uu}=H_{su}=H_{us}=H_{gp}=\{0\}.$ Hence $(V_1, V_2)$ decomposes between pairs given by diagrams and continuously given pairs.
\end{proof}

The continuously given pair do not have such a nice model like pairs given by diagrams or pairs of generalized powers. Hence we do not formulate the result precisely describing the evanescent part for irrational past like we did in previous sections. However, Theorems \ref{Thm_HtHe}, \ref{Thmrs}, \ref{Thmirs} shows that evanescent part of the random field is significantly different in all the three considered pasts.

In our opinion it may be difficult to describe the evanescent part with irrational slope using geometrical construction. It requires rather measure theory approach. It follows from the description of the generator of a continuously given pair. To explain the notion of the generator we need to recall the general idea of the construction of the model of compatible pairs given in \cite{BKS, BKPSCorr}. At first, there were described all the pairs where $\ker V_1^*\cap\ker V_2^*=\{0\},$ which were not the case of a continuously given pair. Then $\ker V_1^*\cap\ker V_2^*$ was decomposed among subspaces generating certain pairs of compatible isometries (including a continuously given pair), where the pair generated by $L\subset \ker V_1^*\cap\ker V_2^*$ is  $(V_1|_{H_L}, V_2|_{H_L})$ for $H_L$  the minimal $V_1, V_2$ reducing subspace containing $L$. The generator of a continuously given pair is the respective subspace $L_0\subset \ker V_1^*\cap\ker V_2^*$. The decomposition of $\ker V_1^*\cap\ker V_2^*$ is obtained by a certain sequence of decompositions of $\ker V_1^*\cap\ker V_2^*$ such that each summand in a given decomposition is a subspace of some summand (precisely one summand) in the previous decomposition. In other words, we build the succeeding decomposition by decomposing some or all summands in previous decomposition. Hence, properly taking one (nontrivial) summand in each decomposition we obtain a decreasing sequence of nontrivial subspaces. Since decompositions are orthogonal, limits of such sequences are pairwise orthogonal as well and only countably many of them may be nontrivial (in separable Hilbert space). The target decomposition of $\ker V_1^*\cap\ker V_2^*$ is the decomposition among nontrivial limits of the described sequences (generating pairs given by diagrams or pairs of generalized powers) and the remaining $L_0$ generating a continuously given pair. The main point is that if there were only countably many decreasing nontrivial sequences we could consider the decomposition among all the limits and then $L_0=\{0\}.$ Hence $L_0\ne\{0\}$ yields an uncountable number of the nontrivial sequences vanishing to $\{0\}$. It reflects the nature of a non-atomic measure. Indeed, if we set a correspondence (possibly constructing some spectral measure) between $L_0$ and some measurable set $A$, such that the  sequence of decompositions of $L_0$  corresponds to the sequence of decompositions of $A$ among countably many disjoint (up to measure $0$) measurable sets (such that each set in a given decomposition is a subset of precisely one set in the previous decomposition), then the fact that each sequence of subspaces vanish to $\{0\}$ corresponds to the fact that each sequence of subsets vanish to a set of measure $0$.

The evanescent part of the random field  with irrational nonsymmetrical half-plane past was investigated by measures in \cite{HL1961, Cuny}. More precisely, in \cite{HL1961} there are characterized purely evanescent processes for irrational nonsymmetrical half-plane past under assumption that the distribution is absolutely continuous with respect to the Lebesgue measure. In \cite{Cuny} there is given an example under similar assumptions, however for the distribution singular to the Lebesgue measure.

\section*{Appendix}

\subsection*{Characterization of half-planes}
In this section we characterize general half-planes. The general idea of a half-plane is to divide points in $\mathbb{Z}^2$ in two parts, by a line. Points on one side of the line are suppose to be in $S,$ while points on the other side in $-S.$ Since $(0,0)$ is neither in $S$ nor in $-S$ it should belong to the line. Moreover, we need the line to pass through $(0,0)$ to get $S\cap-S=\emptyset$. However, there may be also points belonging to the line different from $(0,0)$ which have to be assigned either to $S$ or to $-S$. Let  the line be $\kappa:py-qx=0.$ If $p, q$ are relatively irrational ($p=rq, r\in\mathbb{R}\setminus\mathbb{Q}$), then  $\kappa\cap \mathbb{Z}^2=\{(0,0)\},$ so it is not the case. However, if $p,q$ are relatively rational, assuming $p,q$ to be relatively prime integers, we get  $\kappa\cap \mathbb{Z}^2=\{(kp,kq):k\in\mathbb{Z}\}$. If we assign $(p,q)\in S$, then $S$ as a semi-group contains $\kappa_+=\{(kp,kq):k\in\mathbb{Z}_+\}$ and in turn $-S$ contains $\kappa_-=\{(kp,kq):k\in\mathbb{Z}_-\}.$ If we assign $(p,q)\in -S$ then $\kappa_+, \kappa_-$ are included another way around. Hence, we get  two half-planes defined by the same line (the same vector) $S_{\mathbf{v}}$ and $\widehat{S}_{\mathbf{v}}$ as in \eqref{S_v} and \eqref{S_vhat}. It is clear that $S_{\mathbf{v}}$ and $\widehat{S}_{\mathbf{v}}$ are half-planes and they are different if and only if $p,q$ are relatively rational. Moreover, one can check that $-S_{\mathbf{[p,q]}}=\widehat{S}_{\mathbf{[-p,-q]}}$ and $-\widehat{S}_{\mathbf{[p,q]}}=S_{\mathbf{[-p,-q]}}$ and for $S^X=\{(i,j):(i,-j)\in S\}$ we get $S^X_{\mathbf{[p,q]}}=S_{\mathbf{[p,-q]}}$. In other words $S$ is equal to one of $S_{\mathbf{[p,q]}}$ or $\widehat{S}_{\mathbf{[p,q]}}$ if and only if $-S$ is equal to one of such half-planes if and only if $S^X$ is equal to one of such half-planes.

Next we show that any half-plane is of the form $S_{\mathbf{[p,q]}}$ or $\widehat{S}_{\mathbf{[p,q]}}$. The last condition in Definition \ref{hp} yields that $S$ contains precisely one of sets:$\{(-1,0),(0,-1)\}$, $\{(1,0),(0,-1)\}$, $\{(-1,0),(0,1)\}$, $\{(1,0),(0,1)\}$. We may assume without loss of generality that  $(-1,0),(0,-1)\in S$. Indeed, for the other cases $(-1,0),(0,-1)$ is either in $-S$ or $S^X$ or $-S^X.$ The result for any of such half-planes implies the result for $S$. Since $(-1,0),(0,-1)\in S$, for any $(i,j)\in S$ we get $(i,j)+(\mathbb{Z}_-\cup\{0\})^2\subset S$. Consequently, there are sequences $\{M_i\},\{N_j\}\subset\Zzz\cup\{-\infty,\infty\}$ such that \begin{equation}\label{S}S=\bigcup_{j\in\Zzz}\{(i,j):i\le M_j\}=\bigcup_{i\in\Zzz}\{(i,j):j\le N_i\}\end{equation} as in the picture:
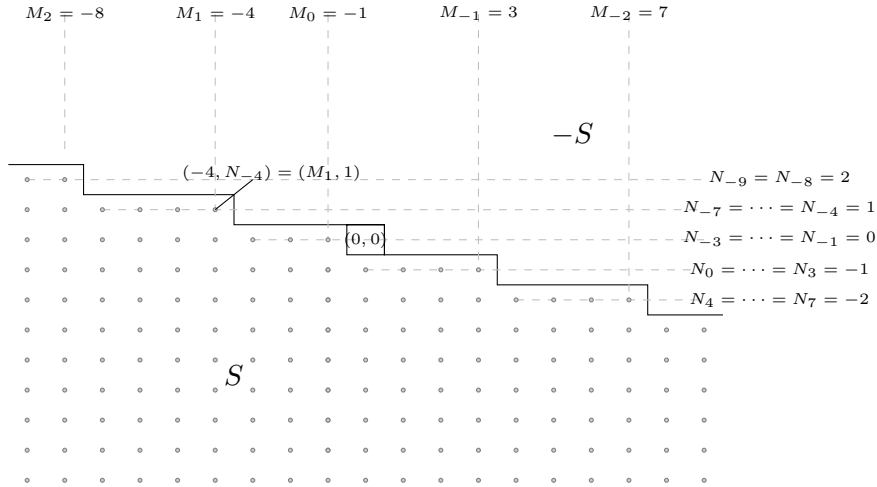
\begin{center}
\begin{tikzpicture}[yscale=0.4, xscale=0.5, auto,
kropka/.style={draw=black!50, circle, inner sep=0pt, minimum size=1.5pt, fill=black!20}]

\draw [black!100] (-10,2) -- ++(2,0) -- ++ (0,-1)-- ++(4,0) -- ++(0,-1) -- ++(4,0) -- ++(0,-1)-- ++(-1,0) -- ++(0,1) -- ++(1,0) -- ++ (0,-1)
 -- ++ (3,0) -- ++ (0,-1) -- ++ (4,0) -- ++ (0,-1)-- ++ (2,0);
\draw [dashed,black!25] (-8.5,7)--(-8.5,2.5);\path(-8.5,7) node[]{\tiny{$M_{2}=-8$}};
\draw [dashed,black!25] (-4.5,7)--(-4.5,0.5);\path(-4.5,7) node[]{\tiny{$M_{1}=-4$}};
\draw [dashed,black!25] (-1.5,7)--(-1.5,-0.5);\path(-1.5,7) node[]{\tiny{$M_{0}=-1$}};
\draw [dashed,black!25] (2.5,7)--(2.5,-1.5);\path(2.5,7) node[]{\tiny{$M_{-1}=3$}};
\draw [dashed,black!25] (6.5,7)--(6.5,-2.5);\path(6.5,7) node[]{\tiny{$M_{-2}=7$}};

\draw [dashed,black!25] (-9.5,1.5)--(8.5,1.5);\path(10.5,1.5) node[]{\tiny{$N_{-9}=N_{-8}=2$}};
\draw [dashed,black!25] (-7.5,0.5)--(8.5,0.5);\path(10.5,0.5) node[]{\tiny{$N_{-7}=\dots=N_{-4}=1$}};
\draw [dashed,black!25] (-3.5,-0.5)--(8.5,-0.5);\path(10.5,-0.5) node[]{\tiny{$N_{-3}=\dots=N_{-1}=0$}};
\draw [dashed,black!25] (-0.5,-1.5)--(8.5,-1.5);\path(10.5,-1.5) node[]{\tiny{$N_{0}=\dots=N_{3}=-1$}};
\draw [dashed,black!25] (3.5,-2.5)--(8.5,-2.5);\path(10.5,-2.5) node[]{\tiny{$N_{4}=\dots=N_{7}=-2$}};

\path (-9.5,1.5) node[kropka]{} ++(0,-1) node [kropka]{}++(0,-1) node [kropka]{}++(0,-1) node [kropka]{}++(0,-1) node [kropka]{}++(0,-1) node [kropka]{}++(0,-1) node [kropka]{}++(0,-1) node [kropka]{}++(0,-1) node [kropka]{}++(0,-1) node [kropka]{}++(0,-1) node [kropka]{};

\path (-8.5,1.5) node[kropka]{} ++(0,-1) node [kropka]{}++(0,-1) node [kropka]{}++(0,-1) node [kropka]{}++(0,-1) node [kropka]{}++(0,-1) node [kropka]{}++(0,-1) node [kropka]{}++(0,-1) node [kropka]{}++(0,-1) node [kropka]{}++(0,-1) node [kropka]{}++(0,-1) node [kropka]{};

\path (-7.5,0.5) node[kropka]{} ++(0,-1) node [kropka]{}++(0,-1) node [kropka]{}++(0,-1) node [kropka]{}++(0,-1) node [kropka]{}++(0,-1) node [kropka]{}++(0,-1) node [kropka]{}++(0,-1) node [kropka]{}++(0,-1) node [kropka]{}++(0,-1) node [kropka]{};

\path (-6.5,0.5) node[kropka]{} ++(0,-1) node [kropka]{}++(0,-1) node [kropka]{}++(0,-1) node [kropka]{}++(0,-1) node [kropka]{}++(0,-1) node [kropka]{}++(0,-1) node [kropka]{}++(0,-1) node [kropka]{}++(0,-1) node [kropka]{}++(0,-1) node [kropka]{};

\path (-5.5,0.5) node[kropka]{}  ++(0,-1) node [kropka]{}++(0,-1) node [kropka]{}++(0,-1) node [kropka]{}++(0,-1) node [kropka]{}++(0,-1) node [kropka]{}++(0,-1) node [kropka]{}++(0,-1) node [kropka]{}++(0,-1) node [kropka]{}++(0,-1) node [kropka]{};

\path (-4.5,0.5) node[kropka]{}  ++(0,-1) node [kropka]{}++(0,-1) node [kropka]{}++(0,-1) node [kropka]{}++(0,-1) node [kropka]{}++(0,-1) node [kropka]{}++(0,-1) node [kropka]{}++(0,-1) node [kropka]{}++(0,-1) node [kropka]{}++(0,-1) node [kropka]{};\draw (-4.5,0.5) -> (-3.5,1.5);
\path (-3,1.7) node[] {\tiny{$(-4,N_{-4})=(M_{1},1)$}};

\path (-3.5,-0.5) node[kropka]{} ++(0,-1) node [kropka] {} ++(0,-1) node [kropka]{}++(0,-1) node [kropka]{}++(0,-1) node [kropka]{}++(0,-1) node [kropka]{}++(0,-1) node [kropka]{}++(0,-1) node [kropka]{}++(0,-1) node [kropka]{};

\path (-2.5,-0.5) node[kropka]{} ++(0,-1) node [kropka] {} ++(0,-1) node [kropka]{}++(0,-1) node [kropka]{}++(0,-1) node [kropka]{}++(0,-1) node [kropka]{}++(0,-1) node [kropka]{}++(0,-1) node [kropka]{}++(0,-1) node [kropka]{};

\path (-1.5,-0.5) node[kropka]{} ++(0,-1) node [kropka] {} ++(0,-1) node [kropka]{}++(0,-1) node [kropka]{}++(0,-1) node [kropka]{}++(0,-1) node [kropka]{}++(0,-1) node [kropka]{}++(0,-1) node [kropka]{};

\path (-1.5,-1.5)  node [kropka] {} ++(0,-1) node [kropka]{}++(0,-1) node [kropka]{}++(0,-1) node [kropka]{}++(0,-1) node [kropka]{}++(0,-1) node [kropka]{}++(0,-1) node [kropka]{}++(0,-1) node [kropka]{};

\path (-0.5,-0.5)  node[]{\tiny{$(0,0)$}} ++(0,-1) node[kropka]{} ++(0,-1) node [kropka] {} ++(0,-1) node [kropka]{}++(0,-1) node [kropka]{}++(0,-1) node [kropka]{}++(0,-1) node [kropka]{}++(0,-1) node [kropka]{}++(0,-1) node [kropka]{};

\path (0.5,-1.5) node[kropka]{} ++(0,-1) node [kropka] {} ++(0,-1) node [kropka]{}++(0,-1) node [kropka]{}++(0,-1) node [kropka]{}++(0,-1) node [kropka]{}++(0,-1) node [kropka]{}++(0,-1) node [kropka]{};

\path (1.5,-1.5) node[kropka]{} ++(0,-1) node [kropka] {} ++(0,-1) node [kropka]{}++(0,-1) node [kropka]{}++(0,-1) node [kropka]{}++(0,-1) node [kropka]{}++(0,-1) node [kropka]{}++(0,-1) node [kropka]{};

\path (2.5,-1.5) node[kropka]{} ++(0,-1) node [kropka] {} ++(0,-1) node [kropka]{}++(0,-1) node [kropka]{}++(0,-1) node [kropka]{}++(0,-1) node [kropka]{}++(0,-1) node [kropka]{}++(0,-1) node [kropka]{};

\path (3.5,-2.5) node[kropka]{} ++(0,-1) node [kropka] {}++(0,-1) node [kropka] {} ++(0,-1) node [kropka]{}++(0,-1) node [kropka]{}++(0,-1) node [kropka]{}++(0,-1) node [kropka]{};

\path (4.5,-2.5) node[kropka]{} ++(0,-1) node [kropka] {}++(0,-1) node [kropka] {} ++(0,-1) node [kropka]{}++(0,-1) node [kropka]{}++(0,-1) node [kropka]{}++(0,-1) node [kropka]{};

\path (5.5,-2.5) node[kropka]{} ++(0,-1) node [kropka] {}++(0,-1) node [kropka] {} ++(0,-1) node [kropka]{}++(0,-1) node [kropka]{}++(0,-1) node [kropka]{}++(0,-1) node [kropka]{};

\path (6.5,-2.5) node[kropka]{} ++(0,-1) node [kropka] {}++(0,-1) node [kropka] {} ++(0,-1) node [kropka]{}++(0,-1) node [kropka]{}++(0,-1) node [kropka]{}++(0,-1) node [kropka]{};

\path (7.5,-3.5) node[kropka]{} ++(0,-1) node [kropka] {}++(0,-1) node [kropka] {}++(0,-1) node [kropka] {} ++(0,-1) node [kropka]{}++(0,-1) node [kropka]{};

\path (8.5,-3.5) node[kropka]{} ++(0,-1) node [kropka] {}++(0,-1) node [kropka] {}++(0,-1) node [kropka] {} ++(0,-1) node [kropka]{}++(0,-1) node [kropka]{};
{\large \path (-4,-5) node{$S$};}
{\large \path (5,3) node{$-S$};}
\end{tikzpicture}
\end{center}
The main idea of the sequences $\{M_i\},\{N_j\}$ is to describe ''corner points'' of $S$ as points $(M_j,j)=(i,N_i)$ for some $i,j$. 

We can see in the picture some regularity. Precisely, both sequences are decreasing and $M_j-M_{j+1}\in\{3,4\}, N_i-N_{i+1}\in\{0,1\}$ for all $i,j\in\mathbb{Z}$. In fact, for any half-plane not being $L$ nor $L'=\{(i,j):(j,i)\in L\}$ there are non-negative integers $m,n$ such that $M_j-M_{j+1}\in\{m,m+1\}$  and $N_i-N_{i+1}\in\{n,n+1\}$ and at least one of $m, n$ is zero. Let us show this fact.

Since $(i,j)+(\mathbb{Z}_-\cup\{0\})^2\subset S$ for any $(i,j)\in S$ we get that $\{M_j\}_{j\in\mathbb{Z}}, \{N_i\}_{i\in\mathbb{Z}}$ are decreasing.  Since $(0,0)\notin S$ and $(-1,0),(0,-1)\in S$ we get $M_0=N_0=-1.$
Moreover,

$$-S=\bigcup_{j\in\Zzz}\{(-i,-j):i\le M_j\}=\bigcup_{j\in\Zzz}\{(i,j):-i\le M_{-j}\}=\bigcup_{j\in\Zzz}\{(i,j):i\ge -M_{-j}\}$$
and
$$-S\cup\{(0,0)\}=\Zzz^2\setminus S=\bigcup_{j\in\Zzz}\{(i,j):i> M_j\}.$$
and similar arguments for $\{N_i\}_{i\in\mathbb{Z}}$ imply \begin{equation}\label{Mi}\begin{array}{ccc}
                                    M_j=-M_{-j}-1&\text{ for }&j\ne 0\\
 N_i=-N_{-i}-1&\text{ for }&i\ne 0.
                                  \end{array}\end{equation}

Let $p\in\Zzz$. Since $(M_{-p},-p), (M_p,p)\in S$ and $M_{-p}=-M_p-1$ by \eqref{Mi}, we get
$$(M_{j+p},j+p)+(M_{-p},-p)=(M_{j+p}-M_p-1,j)\in S\text{ yields }M_{j+p}-M_p-1\le M_j$$ and $$(M_j,j)+(M_p,p)=(M_j+M_p,j+p)\in S,\text{ yields }M_j+M_p\le M_{j+p}$$ for an arbitrary $j$. Summing up \begin{equation}\label{Ni}M_p\le M_{j+p}-M_j\le M_p+1.\end{equation}
In particular for $p=1$ we get $-M_1\ge M_j-M_{j+1}\ge -M_1-1.$ Since $M_1\le M_0=-1$ we get $m=-M_1-1$ a non-negative integer and \begin{equation}\label{m}M_j-M_{j+1}\in\{m,m+1\}.\end{equation}

We may get similarly as above that $N_i-N_{i+1}\in\{n,n+1\}$ for some non-negative integer $n$.

For the further characterization note that one ( and only one ) of elements $(-1,1), (1,-1)$ belongs to $S$. If we assume $(1,-1)\in S$ then $(i, N_i)+(1,-1)=(i+1,N_i-1)\in S$ yields $N_i-1\le N_{i+1}$. Since $\{N_i\}$ is decreasing  $$N_i-1\le N_{i+1}\le N_i.$$ 
In other words, from $(1,-1)\in S$ follows $N_i-N_{i+1}\in\{0,1\}.$ As one can expect the other case $(-1,1)\in S$ yields $M_j-M_{j+1}\in\{0,1\}$.

\begin{proposition}\label{n_k}
For a given $\mathbf{M}=\{M_j\}_{j\in\mathbb{Z}}$ of integers satisfying \eqref{Mi} and \eqref{Ni} for any positive $j, p$ with $M_{0}=-1$ the set 
\begin{equation*}S_{\mathbf{M}}=\left\{(i,j)\in\mathbb{Z}^2:i\le M_j\right\}\end{equation*}
is a half-plane.
\end{proposition}
\begin{proof}
Let us first check that $S_{\mathbf{M}}$ satisfies conditions of Definition \ref{hp}. \begin{itemize}
    \item Since $M_0=-1$ clearly $(0,0)\notin S_{\mathbf{M}}$.
    \item Since by \eqref{Ni} $M_j+M_p\le M_{j+p}$ for any $(x,j), (y,p)\in S_{\mathbf{M}}$ we get $x+y\le M_j+M_p\le M_{j+p}$ so $(x,j)+(y,p)=(x+y,j+p)\in S_{\mathbf{M}}$. Hence $S_{\mathbf{M}}$ is semi-group.
    \item Note that $(i,j)\notin S_{\mathbf{M}}$ if and only if $i>M_j$ which by \eqref{Mi} is equivalent to $-i\le -1-M_j=M_{-j},$ so to $(-i,-j)\in S_{\mathbf{M}}.$
\end{itemize}  
It finishes the proof.
\end{proof} 

\begin{proposition}\label{SL}
Any half-plane $S$ containing $(-1,0),(0,-1)$ is either $L$ (see \eqref{QL}) or $L'=\{(i,j):(j,i)\in L\}$ or \begin{equation}\label{Sdelta}S_{\mathbf{M}}=\left\{(i,j)\in\mathbb{Z}^2:i\le M_j\right\},\end{equation} 
where $M_0=-1$ and $\{M_j\}_{j\in\Zzz}$ satisfy \eqref{Mi} and \eqref{Ni}.
\end{proposition}
\begin{proof}
Since $(-1,0),(0,-1)\in S$ by \eqref{S} we get a sequence $\{M_j\}_{j\in\Zzz}$ such that $S=S_{\mathbf{M}}$ and a sequence $\{N_j\}_{j\in\Zzz}$. 
If $M_1=-\infty$ or $N_1=-\infty$ then $S=L$ or $S=L'$ respectively, otherwise $\{M_j\}_{j\in\mathbb{Z}}, \{N_i\}_{i\in\mathbb{Z}}$ has only finite values by \eqref{m} and similar result for $\{N_i\}_{i\in\mathbb{Z}}$.
\end{proof}

The propositions above 
describes all half-planes as half-planes  of the form $S_{\mathbf{M}}$.
\begin{theorem}
 For any half-plane $S$ there is a vector $\mathbf{v}$ such that $S=S_{\mathbf{v}}$ or $S=\widehat{S}_{\mathbf{v}}$ (see \eqref{S_v}, \eqref{S_vhat}).
\end{theorem}
\begin{proof}
 At the beginning of the paragraph we explained that without lost of generality we may assume $(-1,0),(0,-1)\in S$ . Then by Proposition \ref{SL} it is enough to show the result for $S=S_{\mathbf{M}}$ where $\mathbf{M}=\{M_j\}_{j\in\mathbb{Z}}\subset\mathbb{Z}$. If we denote $\delta_j:=\frac1j(M_0-M_j)=\frac1j(-1-M_j)$ for $j> 0$ ,then $M_j=-j\delta_j-1$ for $j>0$ and by \eqref{Mi} $M_j=-j\delta_{-j}$ for $j<0$. We are going to show that $\delta_j$ is convergent and then the statement for $\mathbf{v}=[-1,-\delta]$ where $\delta=\lim_{j\to\infty}\delta_j.$

From \eqref{m} we get $\delta_j=\frac1j\sum_{k=0}^{j-1}(M_k-M_{k+1})\in[m,m+1].$ Hence $\delta_{sup}:=\limsup_{j\to\infty}\delta_j, \delta_{inf}:=\liminf_{j\to\infty}\delta_j$ are in $[m,m+1]$, so in particular are positive, finite. Denote $\epsilon=\frac13(\delta_{sup}-\delta_{inf})$. Clearly $\delta_j$ is convergent if and only if $\epsilon=0.$

By \eqref{Ni} we get $-1-j\delta_j=M_j\le M_{kj}-M_{(k-1)j}\le M_j+1=-j\delta_j$ which is equivalent to 
$0\le M_{(k-1)j}-M_{kj}-j\delta_j\le 1$. Hence $sj\delta_{sj}-sj\delta_j=\sum_{k=1}^s(M_{(k-1)j}-M_{kj}-j\delta_j)\in[0,s]$ and so $0\le j\delta_{sj}-j\delta_j\le 1,$ in particular $j\delta_j\le j\delta_{sj}.$ 

Let $l>0$ be such that $\delta_l\ge\delta_{sup}-\epsilon.$ Decompose any positive integer $t=sl+t_0$ where $0\le t_0\le l-1$. Recall that $M_j$ is decreasing, by which $j\delta_j$ is increasing. Consequently we get
$$t\delta_t\ge sl\delta_{sl}\ge sl\delta_l\ge (t-t_0)(\delta_{sup}-\epsilon).$$

On the other hand, $t$ may be chosen arbitrary large and  such that $\delta_t\le\delta_{inf}+\epsilon.$ Consequently
$$(t-t_0)(\delta_{sup}-\epsilon)\le t\delta_t\le t(\delta_{inf}+\epsilon)=t(\delta_{sup}-2\epsilon),$$
which yields $$t\epsilon\le t_0(\delta_{sup}-\epsilon)\le (l-1)\delta_{sup}.$$ Since $t$ may be chosen arbitrary large, $\epsilon$ is non negative and the right hand side is constant (for the fixed $l$) we get $\epsilon=0$.

Let us show that $S_{\mathbf{M}}=S_{\mathbf{[-1,-\delta]}}$ or $S_{\mathbf{M}}=\widehat{S}_{\mathbf{[-1,-\delta]}}.$  Since a half-plane may not contain another half-plane unless they are equal it is enough to show   $S_{\mathbf{[-1,-\delta]}}\subset S_{\mathbf{M}}$ or $\widehat{S}_{\mathbf{[-1,-\delta]}}\subset S_{\mathbf{M}}.$ Recall
$$S_{\mathbf{M}}=\left\{(i,j)\in\mathbb{Z}^2:i\le M_j\right\},$$
$$S_{\mathbf{[-1,-\delta]}}=\{(i,j): i\le-j\delta \textnormal{ for } i< 0, i<-j\delta \textnormal{ for } i \ge 0\}.$$

Note that $0\le j\delta_{sj}-j\delta_j\le 1,$ showed earlier yields (taking $s\to\infty$) $j\delta_j\le j\delta\le j\delta_j+1.$ 

Assume $(i,j)\in S_{\mathbf{[-1,-\delta]}}$.
\begin{itemize}
\item if $i\ge 0,$ inequality $i<-j\delta$ may be satisfied only for $j<0$ and we get $i<-j\delta\le-j\delta_{-j}+1=M_j+1$ which yields $i\le M_j$ so $(i,j)\in S_{\mathbf{M}}.$
    \item if $i<0,$  then  $i\le -j\delta\le \left\{
                        \begin{array}{ll}
                          -j\delta_{-j}+1=M_j+1, & \hbox{for $j<0$;} \\
                          -j\delta_j=M_j+1, & \hbox{for $j\ge 0$.}\\
                           \end{array}
                      \right.$ For $j\le 0$ we get $i<-j\delta\leq M_j+1$ which yields $i\le M_j,$ so $(i,j)\in S_{\mathbf{M}}.$ For $j>0$ either $i<M_j+1$ or $i=-j\delta=-j\delta_j=M_j+1.$ In the first case $(i,j)\in S_{\mathbf{M}}$. In the second case, since $i, j$ are integers  $i=-j\delta$ yields $\delta\in\mathbb{Q}.$ Let $\delta=\frac{p}{q}$ for $p, q$ relatively prime integers. Then $i=-j\delta$ yields $(i,j)\in \{(-kp,kq):k\in\mathbb{Z}_+\}.$ In other words, $S_{\mathbf{[-1,-\delta]}}\setminus S_{\mathbf{M}}=S_{\mathbf{[-1,-\delta]}}\cap -S_{\mathbf{M}}\subset \{(-kp,kq):k\in\mathbb{Z}_+\}.$ On the other hand $(-p,q)$ is either in $S_{\mathbf{[-1,-\delta]}}$ or in $-S_{\mathbf{[-1,-\delta]}}$ and either in $S_{\mathbf{M}}$ or in $-S_{\mathbf{M}}.$ Since all of them are semi-groups the whole $\{(-kp,kq):k\in\mathbb{Z}_+\}$ is a subset of those semi-groups where $(-p,q)$ belongs. Hence, either $S_{\mathbf{[-1,-\delta]}}\cap -S_{\mathbf{M}}=\emptyset,$ so $S_{\mathbf{[-1,-\delta]}}\subset S_{\mathbf{M}}$ or $\{(-kp,kq):k\in\mathbb{Z}_+\}\subset S_{\mathbf{[-1,-\delta]}}\cap -S_{\mathbf{M}}.$ 
\end{itemize}

Summing up, either we get $S_{\mathbf{[-1,-\delta]}}= S_{\mathbf{M}},$ so the statement, or $\delta=\frac{p}{q}$ for $p, q$ relatively prime integers and $\{(-kp,kq):k\in\mathbb{Z}_+\}\subset S_{\mathbf{[-1,-\delta]}}\cap -S_{\mathbf{M}}.$
The latter case implies $\{(kp,-kq):k\in\mathbb{Z}_+\}\subset S_{\mathbf{M}}.$ The proof will be done if we show that the last inclusion implies $\widehat{S}_{\mathbf{[-1,-\delta]}}\subset S_{\mathbf{M}}.$ Recall that $$\widehat{S}_{\mathbf{[-1,-\delta]}}=\{(i,j): i<-j\delta \textnormal{ for } i\le 0, i\le-j\delta \textnormal{ for } i > 0\}.$$ For $(i,j)\in \widehat{S}_{\mathbf{[-1,-\delta]}}$:
\begin{itemize}
    \item if $i\le 0,$ then $i< -j\delta\le \left\{
                        \begin{array}{ll}
                          -j\delta_{-j}+1=M_j+1, & \hbox{for $j<0$;} \\
                          -j\delta_j=M_j+1, & \hbox{for $j\ge 0$,}\\
                           \end{array}
                      \right.$ which yields $i\le M_j,$ so $(i,j)\in S_{\mathbf{M}}.$
        \item if $i> 0,$ then  $i\le -j\delta$ is possible only for $j<0$ and $i\le -j\delta\le -j\delta_{-j}+1=M_j+1$ for $j<0$. Hence $i\le M_j+1.$ However, if $i=M_j+1,$ then $i=-j\delta$ and in turn $(i,j)\in\{(kp,-kq):k\in\mathbb{Z}_+\}\subset S_{\mathbf{M}}$ so $i\le M_j$ - a contradiction. Consequently,  $i\leq M_j$, so $(i,j)\in S_{\mathbf{M}}.$
\end{itemize}

\end{proof}

\section*{Data Availability}
Our manuscript has no associated data.

\end{document}